\documentclass[11pt]{article}
\usepackage{latexsym}
\usepackage{amssymb}
\usepackage{bbm}
\usepackage[all]{xy}
\usepackage{graphicx}
\usepackage{amsmath}

\usepackage{diagrams}

\makeatletter
\def\@begintheorem#1#2{\trivlist%
   \item[\hskip \labelsep{\sffamily\bfseries #1\ #2}]\itshape}
\newtheorem{teo}{Theorem}[section]
\newtheorem{defi}[teo]{Definition}
\newtheorem{cor}[teo]{Corollary}
\newtheorem{lem}[teo]{Lemma}
\newtheorem{pro}[teo]{Proposition}
\newtheorem{_rem}[teo]{Remark}
\newtheorem{_eje}[teo]{Example}
\newtheorem{_conj}[teo]{Conjecture}
\newenvironment{rem}{\def\@begintheorem##1##2{\trivlist%
 \item[\hskip\labelsep{\sffamily\bfseries ##1\ ##2}]}\begin{_rem}}{\end{_rem}}
\newenvironment{eje}{\def\@begintheorem##1##2{\trivlist%
 \item[\hskip\labelsep{\sffamily\bfseries ##1\ ##2}]}\begin{_eje}}{\end{_eje}}

\makeatother

\newenvironment{beweis}{{\em Proof:}}{\hfill $\rule{2mm}{2mm}$
\vspace{3mm}

}

\DeclareMathAlphabet{\Ma}{U}{msa}{m}{n}
\DeclareMathAlphabet{\Mb}{U}{msb}{m}{n}
\DeclareMathAlphabet{\Meuf}{U}{euf}{m}{n}
\DeclareSymbolFont{ASMa}{U}{msa}{m}{n}
\DeclareSymbolFont{ASMb}{U}{msb}{m}{n}
\DeclareMathSymbol{\hrist}{\mathord}{ASMa}{"16}
\DeclareMathSymbol{\varkappa}{\mathalpha}{ASMb}{"7B}
\DeclareMathSymbol{\CrPr}{\mathord}{ASMb}{"6F}

\def\got#1{\Meuf{#1}}

\def\mr #1.{\mathrm{#1\,}}
\def\mrt #1.{\mathrm{\mbox{\tiny #1\,}}}
\def\mt #1.{{\mbox{\tiny $#1$}}}
\def\ms #1.{{\mbox{\small $#1$}}}
\def\ol #1.{\overline{#1}}
\def\mb #1.{\mathbf{#1\,}}

\def\1{\mathbbm 1}
\def\wg{\wwh {\al G.}.}

\def\restriction{{\mathchoice{
 \mbox{\unitlength1cm\begin{picture}(.2,.4)%
  \bezier{5}(.07,.3)(.1,.27)(.13,.24)%
  \put(.07,.35){\line(0,-1){.5}}\end{picture}}}{
 \mbox{\unitlength1cm\begin{picture}(.2,.4)%
  \bezier{5}(.07,.3)(.1,.27)(.13,.24)%
  \put(.07,.35){\line(0,-1){.5}}\end{picture}}}{
  \hrist}{\hrist}}}

  \def\al #1.{{\mathcal{#1}}}
  \def\ot #1.{{\got{#1}}}
  \def\C{\Mb{C}}
  \def\c{\mt {\C}.}
  \def\kl #1.{{\mbox{\tiny {$\Mb{#1}$}}}}
  \def\c{\kl C.}

  \def\N{\Mb{N}}
  \def\Z{\Mb{Z}}

  \def\t #1.{\tilde{#1}}
  \def\T #1.{\widetilde{#1}}

\def\bt{\mt {\boxtimes}.}

\def\Wort#1{\mbox{\fontfamily{cmr}\selectfont\mdseries\upshape #1}}

\DeclareMathSymbol{\hsemi}{\mathord}{ASMb}{"6E}
\newcommand{\semi}[2]{\mbox{$#1\kern.1em\hsemi\kern.1em#2$}}
\def\askp#1#2{\mbox{$\LA#1\mbox{\bf ,}\;#2\RA$}}

\def\LA{\left\langle\bgroup}
\def\LE{\left[\bgroup}
\def\LG{\left\{\bgroup}
\def\LR{\left(\bgroup}
\def\RA{\egroup^{\rule{0mm}{2mm}}\right\rangle}
\def\RE{\egroup^{\rule{0mm}{2mm}}\right]}
\def\RG{\egroup^{\rule{0mm}{2mm}}\right\}}
\def\RR{\egroup^{\rule{0mm}{2mm}}\right)}
\def\Ldummy{\left.\bgroup}
\def\Rdummy{\egroup^{\rule{0mm}{2mm}}\right.}

\def\Kbegin{\begin{equation} \left. \begin{array}{rcl}}
\def\Kend{\end{array} \right\} \end{equation}}

\newcommand{\onb}{ortho\-normal basis}
\def\l2{\Lambda^{\mbox{\tiny $(2)$}}}

  \def\ccr #1,#2.{\overline{\Delta(#1,\,#2)}}

  \def\b #1.{{\bf #1}}
  \def\cross#1.{\mathrel{\mathop{\times}\limits_{#1}}}
  \def\C{\Mb{C}}
  \def\N{\Mb{N}}

  \def\Z{\Mb{Z}}
  \def\T{\Mb{T}}

  \def\wh{\widehat}
  
  \def\wwh #1.{\widehat{#1}}
  \def\wt #1.{\widetilde{#1}}

  \def\cross #1.{\mathrel{\raise 3pt\hbox{$\mathop\times\limits_{#1}$}}}
\def\set #1,#2.{\left\{\,#1\;\bigm|\;#2\,\right\}}
\def\b #1.{{\bf #1}}

\def\aut{{\rm Aut}\,}

\def\endo{{\rm End}\,}

\def\Ad{{\rm Ad}\,}
\def\tr{{\rm Tr}\,}

\def\ol #1.{\overline{#1}}
\def\rn#1.{\romannumeral{#1}}

\def\rest{\restriction}
\def\un{\1}

\def\clo{{\rm clo}}
\def\s #1.{_{\smash{\lower2pt\hbox{\mathsurround=0pt $\scriptstyle #1$}}\mathsurround=3pt}}
\def\bra #1,#2.{{\left\langle #1,\,#2\right\rangle_{\al A.}}}

\def\XP#1!{\renewcommand{\baselinestretch}{.7}\marginpar{{\footnotesize #1}\hfil}
\renewcommand{\baselinestretch}{1.5}}
\def\XB{\marginpar{
{\footnotesize\bf Change~starts----}\lower 11pt\hbox{\mathsurround=0pt$
\!\!\displaystyle{
\Bigg\downarrow}$\mathsurround=3pt}}}
\def\XE{\marginpar{{\footnotesize\bf Change~ends-----}\raise 10pt
\hbox{\mathsurround=0pt$
\!\!\displaystyle{
\Bigg\downarrow}$\mathsurround=3pt}}}
\def\HS{{\{\al F.,\,\al G.\}}}

\def\gfin{{\ot G._0}}

\title{\bf Realization of minimal {\it C*}-dynamical systems in terms of
           Cuntz-Pimsner algebras}
\author{
 {\sc Fernando Lled\'o}\footnote{Institute for Pure and Applied Mathematics,
            RWTH-Aachen, Templergraben 55, D-52062 Aachen, Germany (on leave).
 { e-mail: lledo@iram.rwth-aachen.de}}\\[2mm]
 {\footnotesize Department of Mathematics,}      \\[1mm]
 {\footnotesize University Carlos III Madrid}                   \\ [1mm]
 {\footnotesize Avda. de la Universidad 30, E-28911 Legan\'es (Madrid), Spain.}   \\[1mm]
 {\footnotesize flledo@math.uc3m.es}
\and
 {\sc Ezio Vasselli}\\[2mm]
 {\footnotesize Dipartimento di Matematica,}\\[1mm]
 {\footnotesize University of Rome "La Sapienza"}  \\[1mm]
 {\footnotesize P.le Aldo Moro 2, I-00185 Roma, Italy}    \\[1mm]
 {\footnotesize vasselli@mat.uniroma2.it}}

\begin{document}
\maketitle

\begin{center}
{\em Dedicated to Klaus Fredenhagen on his 60th birthday}
\end{center}
\vspace{.3cm}

\begin{abstract}
In the present article we provide several constructions of 
{\it C*}-dynamical systems $(\al F.,\al G.,\beta)$ 
with a compact group $\al G.$ in terms of
Cuntz-Pimsner algebras. These systems have a minimal relative 
commutant of the fixed-point algebra $\al A.:=\al F.^\al G.$ 
in $\al F.$, i.e.
$\al A.'\cap\al F.=\al Z.$, where $\al Z.$ is the center of $\al A.$,
which is assumed to be nontrivial. In addition, we show
in our models that the 
group action $\beta\colon\al G.\to\mr Aut.\al F.$
has full spectrum, i.e. any unitary irreducible 
representation of $\al G.$ is carried by a $\beta_\al G.$-invariant 
Hilbert space within $\al F.$.

First, we give several constructions
of minimal {\it C*}-dynamical systems in terms of
a single Cuntz-Pimsner algebra $\al F.=\al O._\ot H.$ associated to a
suitable $\al Z.$-bimodule $\ot H.$. 
These examples are labeled by the action
of a discrete Abelian group $\ot C.$ (which we call the chain group) 
on $\al Z.$ and by the choice of a suitable class of finite dimensional
representations of $\al G.$. Second, 
we present a more elaborate construction, where now the
{\it C*}-algebra $\al F.$ is generated by a family of Cuntz-Pimsner
algebras. Here the product of the elements
in different algebras is twisted by the chain group
action. We specify the various constructions 
of {\it C*}-dynamical systems for the group $\al G.=\mr SU.\!(N)$,
$N\geq 2$.
\end{abstract}
\begin{flushleft}
{\small
{\sf MSC-classification:} 46L08, 47L80, 22D25 }
\end{flushleft}
\tableofcontents

\section{Introduction}

Duality of groups plays a central role in abstract harmonic analysis.
Its aim is to reconstruct a group $\al G.$ from  
its dual $\widehat{\al G.}$, i.e.~from the set (of equivalence
classes) of continuous, unitary and irreducible representations, endowed with a 
proper algebraic and topological structure. The most famous duality
result is Pontryagin's duality theorem for locally 
compact {\em Abelian} groups.
For compact, not necessarily Abelian, groups there exist also classical
results due to Tannaka and Krein (see \cite{bHewittI,bHewittII}).
Motivated by a long standing problem in quantum field theory,
Doplicher and Roberts came up with a new duality for compact groups
(see \cite{Doplicher89b} as well as M\"uger's appendix
in \cite{pHalvorson06}).
In the proof of the existence of a compact gauge group
of internal gauge symmetries using only a natural set of axioms for the
algebra of observables, they placed the duality of 
compact groups in the framework of C*-algebras. In this situation,
the C*-algebra of local observables $\al A.$ specifies a categorical structure 
generalizing the representation category of a compact group. 
The objects of this category are no longer finite-dimensional Hilbert spaces
(as in the classical results by Tannaka and Krein), but only a certain
semigroup $\al T.$ of unital endomorphisms of the C*-algebra $\al A.$.
In this setting, $\al A.$ has a trivial center,
i.e. $\al Z.:=\al Z.(\al A.)=\C\1$. 
The arrows of the category are the intertwining operators between these
endomorphisms: for any pair of endomorphisms 
$\sigma,\rho\in\al T.$ one defines{\footnote{In this article we will write 
the set arrows ${\mathrm{Hom}}(\rho,\sigma)$
simply by $(\rho,\sigma)$ for each pair $\rho,\sigma$ of objects.}}
\begin{equation}\label{intw}
(\rho,\,\sigma):=\{ X\in\al A.\mid X\rho(A)=\sigma(A)X\;,
\;A\in\al A.\}.
\end{equation}
This category is a natural example of a tensor C*-category, where the norm
of the arrows is the C*-norm in $\al A.$.
The tensor product of objects is defined as composition of endomorphisms
$\rho,\sigma\mapsto \rho\circ\sigma$ 
and for arrows $X_i\in (\rho_i,\sigma_i)$, $i=1,2$, one defines
the tensor product by
\[
 X_1\times X_2:= X_1\rho_1(X_2)\,.
\]
The unit object $\iota$ is the identity endomorphism, which is simple
iff $\al A.$ has a trivial center (since $(\iota,\iota) = \al Z.$).
If $\al A.$ has a trivial center, then the
representation category of $\al G.$ embeds as a full subcategory into the tensor
C*-category of endomorphisms of $\al A.$. The concrete group dual
can be described in terms of an essentially unique C*-dynamical 
system $(\al F.,\al G.,\beta)$, where $\al F.$ is a unital 
C*-algebra containing the original algebra $\al A.$, and the
action of the compact group $\beta\colon\al G.\to\mr Aut.\al F.$
has full spectrum. This means that for any element in the dual
$D\in\widehat{\al G.}$ there is a 
$\beta_\al G.$-invariant Hilbert space $\al H._\mt D.$
in $\al F.$ such that $\beta_\al G.\restriction\al H._\mt D.\in D$. 
(Recall that the scalar product of any pair of elements $\psi,\psi'
\in\al H._\mt D.$ is defined as $\langle\psi,\psi'\rangle:=\psi^*\psi'\in\C\1$
and any orthonormal basis in $\al H._\mt D.$ is a set of 
orthogonal isometries $\{\psi_i\}_i$. The
support of $\al H._\mt D.$ the projection given by
the sum of the end projections, i.e.~$\mr supp.\al H._\mt D.=\sum_i\psi_i\psi_i^*$.) 
Moreover, $\al A.$ is the fixed point algebra of the C*-dynamical system,
i.e.~$\al A.=\al F.^\al G.$ and one has that the 
relative commutant of $\al A.$ in $\al F.$ is minimal, 
i.e.~$\al A.'\cap\al F.=\C\1$. This clearly implies $\al Z.=\C\1$.
The C*-algebra $\al F.$ can also be seen as a crossed product of
$\al A.$ by the semigroup $\al T.$ of endomorphisms of $\al A.$
(cf.~\cite{Doplicher89a}): the endomorphisms $\rho\in\al T.$ (which are inner in
$\al A.$) may be implemented in terms of an orthonormal basis 
$\{\psi_i\}_i\subset\al H.$ in $\al F.$. The endomorphism is 
unital iff the corresponding implementing Hilbert space in $\al F.$ has
support $\1$.  

In a series of articles by Baumg\"artel 
and the first author (cf.~\cite{Lledo97b,Lledo01a,Lledo04a})
the duality of compact groups has been generalized to the case where
$\al A.$ has a nontrivial center, i.e.~$\al Z.\supsetneq\C\1$, and
the relative commutant of $\al A.$ in $\al F.$ remains minimal, i.e.
\begin{equation}
\label{eq_min}
\al A.'\cap\al F.=\al Z.\,.
\end{equation}
(We always have the inclusion $\al Z.\subseteq\al A.'\cap\al F.$.)
We define a {\em Hilbert C*-system} to be 
a {\em C*}-dynamical system
$( \al F. , \al G. , \beta )$ with a 
group action that has full spectrum and for which the Hilbert spaces 
in $\al F.$ carrying the irreducible representations of $\al G.$
have support $\1$
(see Section~\ref{HCSsummary} for a precise definition).
These particular
C*-dynamical systems have a rich structured and many relevant properties
hold, for instance, a Parseval like
identity (cf.~\cite[Section~2]{Lledo04a}).
Moreover, there is an abstract characterization
by means of a suitable {\em non full}
inclusion of {\it C*}-categories $\al T._\c\subset\al T.$, where
$\al T._\c$ is a symmetric tensor category with simple unit, conjugates,
subobjects and direct sums (cf.~\cite{Lledo04a}).
A similar construction appeared in by M\"uger in
\cite{Mueger00}, using crossed products 
of braided tensor *-categories with simple units w.r.t.~a full symmetric 
subcategory.

The C*-dynamical systems $(\al F.,\al G.,\beta)$ in this 
more general context
provide natural examples of tensor C*-categories with a nonsimple unit, since
$(\iota,\iota)=\al Z.$.
The analysis of these kind of categories 
demands the extension of basic notions. For example, a new definition
of irreducible object is needed
(cf.~\cite{Lledo01a,Lledo04a}). In this case the intertwiner
space $(\iota,\iota)\supsetneq\C\1$
is a unital Abelian {\em C*}-algebra and an object $\rho\in\al T.$ is said
to be irreducible if the following condition holds:
\begin{equation}
\label{IrrEnd}
 (\rho,\rho) = 1_\rho \times (\iota,\iota)\,,
\end{equation}
where $1_\rho$ is the unit of the {\it C*}-algebra
$(\rho,\rho)$. In other words, $(\rho,\rho)$ is generated by 
$1_\rho$ as a $(\iota,\iota)$-module.
Another new property that appears in the context of
non-simple units is the action
of a discrete Abelian group on $(\iota,\iota)$. To any irreducible
object $\rho$ one can associate an automorphism $\alpha_\rho\in\mr Aut.\al Z.$
by means of
\begin{equation}
\label{eq_exch}
1_\rho \otimes Z = \alpha_\rho( Z ) \otimes 1_\rho
\ \ , \ \
Z \in \al Z. \ \ .
\end{equation}
Using this 
family of automorphisms $\{\alpha_\rho\}_\rho$ we define
an equivalence relation on $\widehat{\al G.}$,
the dual of the compact group $\al G.$, and the corresponding
equivalence classes become the elements of a 
discrete Abelian group $\ot C.(\al G.)$,
which we call the chain group of $\al G.$.
The chain group is isomorphic to the
character group of the center of $\al G.$ and the map
$\rho\mapsto\alpha_\rho$ induces an action of the chain group on $\al Z.$,
\begin{equation}
\label{eq_ch_ac}
\alpha\colon \ot C.(\al G.) \rightarrow \aut \al Z. \,,
\end{equation}
(see Section~\ref{SecChain}).
The obstruction to have $\al T.$ symmetric is encoded in
the action $\alpha$: 
$\al T.$ is symmetric
if and only if $\alpha$ is trivial (cf.~\cite[Section~7]{Lledo04a}).

These structures are so involved that it is a difficult task to produce 
explicit examples of Hilbert {\em C*}-systems with non-simple unit.
Indeed, up to now it has
been done only for {\em Abelian} groups and in the setting of the
{\it C*}-algebras of the canonical commutation resp.~anticommutation relations
in \cite{Baumgaertel01,Baumgaertel05}. Some {\em indirect} examples based on the
abstract characterization in terms of the inclusion of
{\it C*}-categories $\al T._\c\subset\al T.$, can be found
\cite[Section~6]{Lledo04a}. 

The aim of the present article is to
provide a large class of minimal {\it C*}-dynamical systems and
Hilbert {\em C*}-systems for compact {\it non-Abelian} groups.
These examples are labeled by the action
of the chain group on the unital Abelian {\it C*}-algebra $\al Z.$
given in (\ref{eq_ch_ac}). A crucial part of our examples are the
Cuntz-Pimsner algebras introduced by Pimsner in 
his seminal article \cite{Pimsner97}. This is 
a new family of {\it C*}-algebras
$\al O._\al M.$ that are naturally generated by a Hilbert bimodule
$\al M.$ over a {\it C*}-algebra $\al A.$. These algebras
generalize Cuntz-Krieger algebras as well as crossed-products by the group
$\Z$. In Pimsner's construction $\al O._\al M.$ is given as
a quotient of a Toeplitz like algebra acting on a concrete
Fock space associated to $\al M.$.
An alternative abstract approach to Cuntz-Pimsner algebras
in terms of {\it C*}-categories is given in
\cite{Doplicher98,Kajiwara98,PinzariIn97}.
In our models we Cuntz-Pimsner algebras
$\al O._\ot H.$ associated to a certain free $\al Z.$-bimodules
$\ot H.=\al H.  \otimes  \al Z.$. The factor $\al H.$ denotes a generating finite 
dimensional Hilbert space with an orthonormal basis specified by isometries 
$\{\psi_i\}_i$. 
The left $\al Z.$-action of the bimodule is defined in terms of
the chain group action (\ref{eq_ch_ac}).

\subsection{Main results}
To state our first main result we need to introduce the family $\gfin$
of all finite-dimensional representations $V$ of the compact group
$\al G.$ that satisfy the following two properties: first,
$V$ admits an irreducible subrepresentation of dimension
or multiplicity $\geq 2$ and, second,
there is a natural number 
$n \in \N$ such that $\mathop{\otimes}\limits^n V$ contains
the trivial representation $\iota$,
i.e.~$\iota \prec \mathop{\otimes}\limits^n V$.
Then we show:\\[3mm]
{\bf Main Theorem~1} (Theorem~\ref{thm_main})
{\it Let $\al G.$ be a compact
group, $\al Z.$ a unital Abelian C*-algebra
and consider a fixed chain group action
$\alpha  \colon \ \ot C.(\al G.) \to \mr Aut.(\al Z.)$. Then for
any $V\in\gfin$ there exists a $\al Z.$-bimodule $\ot H._\mt V.
=\al H._\mt V.  \otimes  \al Z.$ 
with left $\al Z.$-action given in terms
of $\alpha$ 
and a C*-dynamical system $(\al O._{\ot H._\mt V.} ,\al G.,\beta_\mt V.)$,
satisfying the following properties:
\begin{itemize}
\item[(i)] $(\al O._{\ot H._\mt V.} ,\al G.,\beta_\mt V.)$ is minimal,
i.e.~$\al A._\mt V.' \cap \al O._{\ot H._\mt V.} = \al Z.$, where
$\al A._\mt V. := \al O._{\ot H._\mt V.}^{\al G.}$
is the corresponding fixed-point algebra.
\item[(ii)] The Abelian C*-algebra $\al Z.$ coincides with
the center of the fixed-point algebra $\al A._\mt V.$, i.e.~
$\al A._\mt V.' \cap \al A._\mt V. = \al Z.$.
\end{itemize}
Moreover, if $\al G.$ is a compact Lie group, then the Hilbert
spectrum of $(\al O._{\ot H._\mt V.}  ,\al G.,\beta_\mt V.)$ is full, i.e.~for each
irreducible class $D\in\wh{\al G.}$ there is an invariant Hilbert space
$\al H._{\mt D.} \subset\al O._{\ot H._\mt V.}$
(in this case not necessarily of support $\1$)
such that $\beta_\mt V.\restriction \al H._{\mt D.}$ specifies an
irreducible representation of class $D$.}\\[3mm]
An important step in the proof is to show that the corresponding bimodules
$\ot H._\mt V.$ are nonsingular. This notion was introduced in
\cite{Doplicher98} and is important for analyzing the relative commutants in
the corresponding Cuntz-Pimsner algebras (see Section~\ref{Cuntz_Pimsner} for
further details). We give a characterization of the class of nonsingular 
bimodules that will appear in this article (cf.~Proposition~\ref{lem_nsing}).
The preceding theorem may be applied to the group
SU$(N)$ in order to define a corresponding minimal {\it C*}-dynamical system
with full spectrum (cf.~Example~\ref{SUN}).

To present examples of minimal {\it C*}-dynamical systems 
with full spectrum, where the Hilbert spaces in $\al F.$ that carry the 
irreducible representations of the group have support $\1$, we need
a more elaborate construction: to begin with,
we introduce a {\it C*}-algebra generated by a family of Cuntz-Pimsner
algebras that are labeled by {\em any} family $\ot G.$ of
unitary, finite-dimensional representations of $\al G.$
(see Subsection~\ref{5.1} for precise presentation
of this algebra). This construction is interesting in itself and can
be performed for coefficient algebras $\ot R.$ which are not necessarily
Abelian. Concretely we show:\\[3mm]
{\bf Main Theorem~2} (Theorem~\ref{thm_exf}) {\it
Let $\al G.$ be a compact group, $\ot R.$ a unital {\it C*}-algebra
and $\alpha\colon\ot C.(\al G.) \rightarrow\aut \ot R.$ a fixed
action of the chain group $\ot C.(\al G.)$.
Then, for every set $\ot G.$ of finite-dimensional representations of $\al G.$,
there exists a universal {\it C*}-algebra
$\ot R. \rtimes^\alpha \ot G.$
generated by $\ot R.$ and the Cuntz-Pimsner algebras
$\left\{\al O._{\ot H._\mt V.}\right\}_{\mt {V\in \ot G.}.}$,
where the product of the elements in the different algebras
is twisted by the chain group action $\alpha$.
}\\[3mm]
The {\it C*}-algebra $\ot R. \rtimes^\alpha \ot G.$
(which we will also denote simply by $\al F.$) generalizes some
well-known constructions, obtained for particular choices of the family of
representations $\ot G.$, such as Cuntz-Pimsner algebras, crossed products
by single endomorphisms ({\it \`a la} Stacey) or crossed products by Abelian
groups. Hilbert space representations of $\ot R. \rtimes^\alpha \ot G.$
are labeled by covariant representations of the {\it C*}-dynamical
system $( \ot R. , \ot C. (\al G.) , \alpha )$.

Now, we restrict the result of the Main Theorem~2
to the case $\ot G.=\gfin$
with Abelian coefficient algebra $\ot R.=\al Z.$. 
The C*-algebra $\al F.=\al Z. \rtimes^\alpha \gfin$ specifies
prototypes of Hilbert {\it C*}-systems for non-Abelian groups
in the context of non-simple units satisfying all the 
required properties:\\[3mm]
{\bf Main Theorem~3} (Theorem~\ref{FinalEx}) {\it
Let $\al G.$ be a compact group, $\al Z.$ a unital Abelian {\it C*}-algebra
and $\alpha\colon\ot C.(\al G.) \rightarrow\aut \al Z.$ a fixed
chain group action. Given the set of finite-dimensional
representations $\gfin$ introduced above
and the {\it C*}-algebra $\al F.:=\al Z.\rtimes^\alpha\gfin$
of the preceding theorem, there exists a minimal
{\it C*}-dynamical system $( \al F. , \al G. , \beta )$,
i.e.~$\al A.'\cap\al F.=\al Z.$, where $\al A.$ is the corresponding fixed
point algebra.
Moreover, $\al Z.$ coincides with the center of $\al A.$, i.e.
$\al Z.=\al A.'\cap\al A.$, and
for any $V \in \gfin$ the Hilbert space $\al H._\mt V.
\subset\al O._{\ot H._\mt V.}\subset\al F.$ has support $\1$.
}\\[3mm]
We may apply the preceding theorem to the group $\al G.:=$SU(2).
Here we choose as
the family of finite-dimensional representations $\ot G._0$ all irreducible
representations of $\al G.$ with dimension $\geq 2$. This gives an explicit
example of a Hilbert {\it C*}-system for SU(2) (cf.~Example~\ref{SU2}).

The structure of the article is as follows: In Section~\ref{Sec-2}
we present the main definitions and results concerning
Hilbert {\it C*}-systems and the chain group.
In Section~\ref{Cuntz_Pimsner} we recall the main features of
Cuntz-Pimsner algebras that will be needed later. In the following
section we present a family of minimal {\it C*}-dynamical
systems for a compact group $\al G.$ and a single Cuntz-Pimsner algebra.
This family of examples is labeled by the chain
group action (\ref{eq_ch_ac}) and the elements of a
suitable class $\ot G._0$ of finite-dimensional representations
of $\al G.$.
In Section~\ref{SecCompact} we construct first a {\it C*}-algebra
$\al F.$ generated by the Cuntz-Pimsner algebras
$\{\al O._{\ot H._V}\}_{V\in\ot G._0}$ as described
above. Then we show that with $\al F.$ we
can construct a Hilbert {\it C*}-system in a natural way.
We conclude this article with a short appendix restating some
of the previous concrete results in terms of
tensor categories of Hilbert bimodules.

\subsection{Outlook}

Doplicher and Roberts show in the setting of the new duality of compact groups
that essentially every concrete dual of a compact group
$\al G.$ may be realized in a natural way within a C*-algebra $\al F.$,
which is the C*-tensor product of Cuntz algebras (cf.~\cite{Doplicher88}).
Under additional assumptions it is shown that the corresponding fixed point
algebra is simple and therefore must have a trivial center.
The results in this paper generalize this situation. In fact, 
one may also realize concrete group duals within the 
C*-algebra $\al F.:=\al Z.\rtimes^\alpha\gfin$ constructed in 
the Main Theorem~3, where now
the corresponding fixed point algebra has a nontrivial center $\al Z.$.
If $\al Z.=\C\1$, then $\al Z.\rtimes^\alpha\ot G.$ reduces to the 
tensor product of Cuntz algebras labeled by the finite dimensional
representations of the compact group contained in $\ot G.$.

As mentioned above our models provide natural examples of
tensor C*-categories with a nonsimple unit. These structures
have been studied recently 
in several problems in mathematics and mathematical physics:
in the general context of 2-categories (see
\cite{Zito06} an references cited therein), in the study of group
duality and vector bundles \cite{Vasselli05,Vasselli06},
and in the context of superselection theory
in the presence of quantum constraints \cite{Baumgaertel05}.
Finally, algebras
of quantum observables with nontrivial center 
$\al Z.$ also appear
in lower dimensional quantum field theories with braiding symmetry
(see e.g.~\cite{Fredenhagen92}, \cite[\S 2]{Mack90}).
In particular, in the latter reference the vacuum representation of the 
global observable algebra is not faithful and maps central elements 
to scalars. In the mathematical setting of this article,
the analogue of the observable algebra is analyzed 
without making use of Hilbert space representations that
trivialize the center. Moreover, the representation theory of a
compact group is described by endomorphisms (i.e. the analogue
of superselection sectors) that preserve the center. 
It is clear that our models do not fit completely in
the frame given by lower dimensional quantum field theories, since,
for example, we do not use any braiding symmetry.
Nevertheless, we hope that some pieces of the analysis
considered here can also be applied.
E.g.~the generalization of the notion of irreducible 
objects and the analysis of their restriction to
the center $\al Z.$ that in our context led to the 
definition of the chain group or the importance of Cuntz-Pimsner algebras
associated to $\al Z.$-bimodules.

\section{Hilbert {\it C*}-systems and the chain group}
\label{Sec-2}

For convenience of the reader we recall the main definitions and
results concerning Hilbert {\it C*}-systems
that will be used later in the construction of the examples.
We will also introduce the notion of the chain group associated
to a compact group which will be crucial in the specification of
the examples.
For a more detailed analysis of Hilbert {\it C*}-systems
we refer to \cite[Sections~2 and 3]{Lledo04a} and
\cite[Chapter~10]{bBaumgaertel92}).

\subsection{Hilbert {\it C*}-systems}\label{HCSsummary}

Roughly speaking, a Hilbert {\it C*}-system is a special
type of {\it C*}-dynamical system $\{\al F.,{\al G.},\beta\}$ that,
in addition, contains the information of the representation category
of the compact group $\al G.$.
$\al F.$ denotes a unital {\it C*}-algebra and
$\beta\colon\al G.\ni g \mapsto \beta_g \in \mr Aut.\al F.$ is a pointwise
norm-continuous morphism. Moreover, the representations
of $\al G.$ are carried by the algebraic Hilbert spaces,
i.e.~Hilbert spaces $\al H.\subset\al F.$, where the
scalar product $\langle\cdot,\cdot\rangle$ of $\al H.$
is given by
$\langle A,B\rangle\1 := A^{\ast}B$ for $A,\; B\in\al H.$.
(Algebraic Hilbert spaces are also called in the literature 
Hilbert spaces in C*-algebras.)
Henceforth, we consider only finite-dimensional algebraic Hilbert spaces.
The support $\hbox{supp}\,\al H.$
of $\al H.$ is defined by\
$\hbox{supp}\,\al H.:=\sum_{j=1}^{d}\psi_j\psi_{j}^{\ast}$,
where $\{\psi_j\,\mid j=1,\ldots,\,d\}$ is any orthonormal
basis of $\al H.$.

To give a precise definition of a Hilbert {\it C*}-system we need to
introduce the spectral projections: for $D\in\wh{\al G.}$
(the dual of $\al G.$)
its {\it spectral projection} $\Pi_{\mt D.}\in \al L. (\al F.)$
is defined
by
\begin{eqnarray} \label{PiD}
\Pi_{{\mt D.}} (F)&:=&\int_{\al G.}\ol\chi_{\mt D.} (g).\,\beta_{g}(F)\,dg
\quad\hbox{for all}\quad F\in\al F.,  \\[1mm]
\nonumber \hbox{where}\quad\qquad
\chi_{\mt D.} (g)&:=&\dim{D}\cdot\tr {U_\mt D.}(g),\quad
                    U_\mt D.\in D\,,
\end{eqnarray}
is the so-called modified character of the class $D$ and
$dg$ is the normalized Haar measure of the compact group $\al G.$.
For the trivial representation $\iota\in\widehat{{\cal G}}$, we put
\[
 {\cal A}:=\Pi_{\iota}\,{\cal F}=\LG F\in{\cal F}\mid
             g(F)=F,\quad g\in {\cal G}\RG ,
\]
i.e.~${\cal A}=\al F.^\al G.$ is
the fixed-point algebra in ${\cal F}$
w.r.t.~${\cal G}$. We denote by $\al Z.=\al Z.(\al A.)$ the center of
$\al A.$, which we assume to be nontrivial.

\begin{defi}\label{defs2-1}
The {\it C*}-dynamical
system $\{\al F.,\al G.,\beta\}$ with compact group $\al G.$
is called a {\bf Hilbert {\it C*}-system} if
it has full Hilbert spectrum, i.e.~for each $D\in\wh{\al G.}$
there is a $\beta$-stable Hilbert space
$\al H._{\mt D.}\subset\Pi_{\mt D.}\al F.,$
with support $\1$
and the unitary representation
$\beta_\al G.\rest\al H._{\mt D.}$
is in the equivalence class
$D\in\wh{\al G.}$. A Hilbert {\it C*}-system is
called {\bf minimal} if
\[
 \al A.'\cap\al F.=\al Z. \,,
\]
where $\al Z.$ is the center of the fixed-point algebra
$\al A.:=\al F.^\al G.$.
\end{defi}

Since we can identify $\al G.$ with
$\beta_\al G.\subseteq\mr Aut.\al F.$ we will often denote
the Hilbert {\it C*}-system simply by $\{\al F.,\al G.\}$.

\begin{rem}\label{ExCAR}
Some families of examples of minimal Hilbert {\it C*}-systems with fixed point
algebra $\al A._\c\otimes\al Z.$, where $\al A._\c$ has trivial
center, were constructed indirectly in \cite[Section~6]{Lledo04a}.
Some explicit examples
in the context of the CAR/CCR-algebra
with an Abelian group are given in \cite{Baumgaertel01}
and \cite[Section~V]{Baumgaertel05}.
\end{rem}

To each $\al G.$-invariant algebraic Hilbert space $\al H.\subset
\al F.$
there is assigned a corresponding
inner endomorphism $\rho\s{\al H.}.\in\endo\al F.$ given by
\[
\rho_{\mt {\al H.}.}(F):=\sum_{j=1}^{d(\al H.)}\psi_jF\psi_j^*\;,
\]
where $\{\psi_j\,\big|\,j=1,\ldots,\,d(\al H.)\}$ is any orthonormal basis of
$\al H.$. It is easy to see that $\al A.$ is stable under the inner
endomorphism $\rho$.
We call {\bf canonical endomorphism} the restriction of
$\rho_{\al H.}$ to
$\al A.$,
i.e.
$\rho_{\al H.}\rest\al A.\in\mbox{End}\,\al A.$. By abuse of
notation we will also denote
it simply by $\rho_{\al H.}$. 
Let $\al Z.$ denote the center of $\al A.$;
we say that an endomorphism $\rho$ is {\em irreducible} if
\[
 (\rho,\rho)=\rho(\al Z.)\,.
\]
In the nontrivial center situation 
canonical endomorphisms do not characterize the
algebraic Hilbert spaces anymore. In fact,
the natural generalization in this context is the following notion
of free Hilbert $\al Z.$-bimodule:
let ${\cal H}$ be a ${\cal G}$-invariant algebraic Hilbert space in
${\cal F}$ of finite dimension $d$. Then we define first the free right
${\cal Z}$-module ${\got H}$ by extension
\begin{equation}\label{GenBiMod}
 {\got H}:= {\cal H}{\cal Z}=\LG\sum\limits_{i=1}^d \psi_i\, Z_i\mid
            Z_i\in {\cal Z}\RG,
\end{equation}
where $\Psi:=\{ \psi_i \}_{i=1}^d$ is an \onb{} in ${\cal H}$. In other words,
the set $\Psi$ becomes a module basis of ${\got H}$ and
$\Wort{dim}_{\mt {\al Z.}.}{\got H}=d$. For $H_1,H_2\in{\got H}$ put
\[
 \askp{H_1}{H_2}_{{\got H}}:=H_1^*H_2\in{\cal Z}\,.
\]
Then, $\{ {\got H},\,\askp{\cdot}{\cdot}_{{\got H}}\}$ is a Hilbert right
${\cal Z}$-module or a Hilbert $\al Z.$-module, for short.
Now the canonical endomorphism can be also written as
\[
\rho_{\mt {\al H.}.} (A):=\sum_{j=1}^{d}\varphi_j A \varphi_j^*
\;,\quad A\in\al A.\,,
\]
where $\{ \varphi_i \}_{i=1}^d$ is {\em any} orthonormal basis
of the $\al Z.$-module $\ot H.$. Hence we actually have
$\rho_{\mt {\al H.}.}=\rho_{\mt {\ot H.}.}$ and it is
easy to show that
\[
 H\in\ot H.\quad\mr iff.\quad H\,A=\rho_{\mt {\ot H.}.}(A)\,H\,.
\]
In other words $\rho_{\mt {\ot H.}.}$ characterizes uniquely
the Hilbert $\al Z.$-module $\ot H.$.
Moreover, since for any canonical endomorphism
$\rho=\rho_{\mt {\al H.}.}$
we have that $\al Z.\subset (\rho,\rho)$, it is easy to see that
there is a canonical left action of $\al Z.$ on $\ot H.$.
Concretely, there is a natural *-homomorphism
$\al Z.\to \al L.(\ot H.)$, where $\al L.(\ot H.)$ is the
set of $\al Z.$-module morphisms
(see \cite[Sections~3 and 4]{Lledo97b} for more details).
Hence $\ot H.$ becomes a $\al Z.$-bimodule.

We conclude stating the isomorphism between the category of
canonical endomorphisms and the corresponding
category of free $\al Z.$-bimodules
(cf.~\cite[Proposition~4.4]{Lledo04a} and \cite[Section~4]{Lledo97b}).

\begin{pro}\label{prop0}
Let $\HS$ be a given minimal Hilbert C*-system, where the fixed
point algebra $\al A.$ has center $\al Z.$.
Then the category $\al T.$ of all canonical endomorphisms of
$\HS$ is isomorphic to the subcategory $\al M._\al G.$
of the category of free Hilbert
$\al Z.$-bimodules with objects $\ot H.=\al H.\al Z.$,
where $\al H.$ is a $\al G.$-invariant algebraic
Hilbert space with $\mathrm{supp}\al H._{\sigma}=\1$, and
the arrows given by the corresponding $\al G.$-invariant
module morphisms
$\al L.({\got H}_{1},{\got H}_{2};\al G.)$.

The bijection of objects is given by
$\rho_\al H.\leftrightarrow \ot H.=\al H.\al Z.$
which satisfies the conditions
\begin{eqnarray*}
\rho_\al H.=(\Ad V)\circ\rho_1+(\Ad W)\circ\rho_2
&\longleftrightarrow&
{\got H}=V{\got H}_1+W{\got H}_2\\[1mm]
\rho_1\circ\rho_2 &\longleftrightarrow&
{\got H}_1\cdot{\got H}_2    \;,
\end{eqnarray*}
where $V,W\in\al A.$ are isometries with $VV^*+WW^*=\1$ 
and the latter product is the inner tensor product of the
Hilbert $\al Z.$-modules w.r.t.~the *-homomorphism
$\al Z.\to\al L.(\ot H._2)$.
The bijection on arrows is defined by
\[
\al J.\colon\ \al L.({\got H}_{1},{\got H}_{2};\al G.)
\to (\rho_1,\rho_2) \quad\mr with.\quad
\al J.(T):=\sum_{j,k}\psi_{j}Z_{j,k}\varphi_{k}^{\ast}\;.
\]
Here
$\{\psi_{j}\}_{j},\;\{\varphi_{k}\}_{k}$
are orthonormal basis of
${\got H}_{2},{\got H}_{1}$,
respectively, and
$(Z_{j,k})_{j,k}$
is the matrix of the right $\al Z.$-linear operator $T$ from
${\got H}_{1}$
to
${\got H}_{2}$ which intertwines the $\al G.$-actions.
\end{pro}

The preceding proposition shows that the canonical endomorphisms
uniquely determine the corresponding $\al Z.$-bimodules, but not the choice
of the generating algebraic Hilbert spaces. The assumption
of the minimality condition in Definition~\ref{defs2-1} is crucial here.
From the point of view of the $\al Z.$-bimodules it is
natural to consider next the following property of Hilbert
C*-systems:
the existence of a special choice of algebraic Hilbert spaces
within the modules that define
the canonical endomorphisms and which is compatible with products.

\begin{defi}\label{ReguCondi}
A Hilbert C*-system
$\HS$
is called {\bf regular} if there is an assignment
$\al T.\ni\sigma\rightarrow\al H._{\sigma}$,
where $\al H._{\sigma}$ is a $\al G.$-invariant algebraic Hilbert space
with $\mathrm{supp}\al H._{\sigma}=\1$ and
$\sigma=\rho\s{\al H._\sigma}.$
(i.e.~$\sigma$ is the canonical endomorphism of
the algebraic Hilbert space $\al H._\sigma$),
which is compatible with products:
\[
\sigma\circ\tau \mapsto \al H._\sigma\cdot\al H._\tau\;.
\]
\end{defi}

\begin{rem}
In a minimal Hilbert C*-system regularity means that there is a
``generating" Hilbert space
$\al H._{\tau}\subset{\got H}_{\tau}$
for each $\tau$ (with
${\got H}_{\tau}=\al H._{\tau}\al Z.$)
such that the compatibility relation for products
stated in Definition~\ref{ReguCondi} holds.
If a Hilbert C*-system is minimal and
$\al Z.=\C\un$
then it is necessarily regular.
\end{rem}

\subsection{The chain group}\label{SecChain}

In the present section we recall the main motivations and definitions
concerning the chain group associated with a compact group $\al G.$.
For proofs and more details see \cite[Section~5]{Lledo04a}
(see also \cite{Mueger04}).

One of the fundamental new aspects of superselection theory with
a nontrivial center $\al Z.$
is the fact that irreducible canonical endomorphisms act
as (nontrivial) automorphisms on $\al Z.$.
In fact, let $D\in\wh{\al G.}$ (the dual of $\al G.$)
and denote by $\rho_\mt D.:=\rho_{\al H._D}$
the corresponding irreducible canonical endomorphism. Then, to
any class $D$ we can associate
the following automorphism on $\al Z.$:
\begin{equation}\label{alphaD}
 \wh{\al G.}\ni D \mapsto \alpha_\mt D.:=\rho_\mt D.\rest\al Z.\in\mr Aut.\al Z.\,.
\end{equation}
This observation allows one to introduce a natural
equivalence relation in the dual $\wwh {\al G.}.$
which, roughly speaking, relates elements $D,D'\in\wwh {\al G.}.$ if there
is a ``chain of tensor products'' of elements in $\wwh {\al G.}.$
containing $D$ and $D'$
(see Theorem~\ref{EndoChain} and Remark~\ref{CompletePicture}
below).

First, we need to recall the algebraic structure of $\wwh {\al G.}.$:
denote by $\times$ the natural operation
on subsets of $\wwh {\al G.}.$ associated with the
decomposition of the tensor products of irreducible representations.
For any $D\in\wwh {\al G.}.$ let $U_\mt D.$ be an irreducible representation 
in the class $D$. Then we define
\[
 D_1\times D_2:=\{D\in\wwh {\al G.}.\mid
                  U_\mt D.\prec U_\mt D_1. \otimes U_\mt D_2. \}\;.
\]
For $\Gamma,\Gamma_1,\Gamma_2\subset\wg$ put
\[
 \Gamma_1\times \Gamma_2=\cup\{D_1\times D_2\mid D_i\in \Gamma_i,\,i=1,2\}
 \quad\mr and.\quad D\times\Gamma =\{D\}\times \Gamma\,.
\]
Moreover if $\ol D.\in\wg$ denotes the conjugate class to
$D\in\wg$ we put $\ol \Gamma.=\{\ol D.\mid D\in\Gamma\}$.
Recall in particular that if $D\in D_0\times D_1$,
$D'\in D_0'\times D_1'$, then
$D\times D'\subset D_0\times D_1\times D_0'\times D_1'$
or that the trivial representation $\iota$ is contained in
$\Gamma\times\ol \Gamma.$
(cf.~\cite[Definition~27.35]{bHewittII} for further details).

We can now make precise the previous idea:
\begin{defi}\label{DefChain} The elements
$D,D'\in\wwh {\al G.}.$ are called equivalent, $D\approx D'$,
if there exist $D_1,\ldots, D_n \in\wwh {\al G.}.$ such that
\[
D,D'\in D_1\times\dots\times D_n\,.
\]
\end{defi}
The preceding definition is an equivalence relation in $\wwh {\al G.}.$
and we denote by square brackets $[\cdot]$
the corresponding chain equivalence classes.
We denote the factor space by
\[
  \ot C.(\al G.):= \wwh {\al G.}. /\! \approx\;.
\]
By definition {\em any} pair
$D,D'\in D_0\times D_1$ satisfies
$D\approx D'$. Therefore for
$D_0,D_1\in\wwh {\al G.}.$ we have that $D_0\times D_1$
also specifies an element of $\ot C.(\al G.)$ and
we can simply put
\[
 [D_0\times D_1]:=[D]\,,
\]
where $D$ is any element in $D_0\times D_1$.

We will define on $\ot C.(\al G.)$ a product $\bt$
(see Eq.~(\ref{ChainProduct}) below) so that
$(\ot C.(\al G.),\bt)$ becomes an Abelian group which for simplicity
we call {\bf chain group}.
Moreover, the chain group
can be related to the character group of the
center $\al C.$ of $\al G.$.
For this recall also the notion of
conjugacy class of a representation
(cf.~\cite{bFuchs97}):
let $D\in\wwh {\al G.}.$ and $U_\mt D.$ any representer in
$D$. By Schur's Lemma we
have
\begin{equation}\label{RepCC}
 U_\mt D.\rest\al C.=\Upsilon_\mt D.\cdot\1\;,
\end{equation}
and it can be easily seen that $\Upsilon_\mt D.$ is a character
on the center $\al C.$ of $\al G.$
which only depends on $D$, i.e.~$\Upsilon_\mt D.\in\wwh {\al C.}.$.

\begin{teo}\label{TeoChain}
Let $\al G.$ be a compact non-Abelian group and denote
its center by $\al C.$.
  \begin{itemize}
  \item[(i)] The set $\ot C.(\al G.)$ becomes an Abelian group w.r.t.~the
             following multiplication: for $D_0,D_1\in\wwh {\al G.}.$ put
\begin{equation}\label{ChainProduct}
 [D_0]\;\bt\; [D_1]:= [D_0\times D_1]\,.
\end{equation}
  \item[(ii)] The conjugacy classes $\Upsilon_\mt D.$
              (cf.~Eq.~(\ref{RepCC})) depend
              on the chain equivalence class $[D]$.
              The chain group and the character group of the center of
              $\al G.$ are isomorphic. The isomorphism is given by
\[
 \eta:\ot C.(\al G.)\to\wh{\al C.}
 \quad\mr with.\quad \eta([D]):=\Upsilon_\mt [D].\;,
\]
where $\Upsilon_\mt [D].$ is the conjugacy class associated with
$[D]\in\ot C.(\al G.)$.
\end{itemize}
\end{teo}

\begin{rem}
\begin{itemize}
\item[(i)] A complete proof of the preceding theorem is given in
Theorem~5.5 of \cite{Lledo04a}. The injectivity of the mapping $\eta$ 
was proven for the first time in \cite{Mueger04}.
\item[(ii)]
The preceding theorem shows that the equivalence relation defined in $\wg$
gives a direct way to reconstruct (via Pontryagin's duality)
the center $\al C.(\al G.)$ of a compact group from the representation ring of
$\al G.$.
\begin{diagram}
 \al G.     & \rTo & \al C.(\al G.) \\
 \dTo       &      & \uTo_{\mbox{\tiny Pontryagin}}           \\
 \wg        & \rTo_{\mbox{\tiny factorization by $\approx$}} & \ot C.(\al G.)
\end{diagram}
For further results in this direction see \cite{Mueger04,ZimborasIn06}.
\end{itemize}
\end{rem}

\begin{eje}\label{SU2_v1}
A simple example that illustrates how to construct a chain group
is given by $\al G.=\mr SU.\!(2)$. Denote by
\[
 l\in\{0,\ms \frac12.,1,\ms \frac32.,\dots\}=\wwh {\mr SU.\!(2)}.
\]
the class specified by the usual
representation $V^{(l)}$ of $\mr SU.\!(2)$ on the space of complex polynomials
of degree $\leq 2l$ which has dimension $2l+1$. Then the
decomposition theory for the tensor products $V^{(l)}\otimes V^{(l')}$
(cf.~\cite[Theorem~29.26]{bHewittII}) gives
\[
 l\times l'=\left\{\, |l-l'|,|l-l'|+1,\dots,l+l' \,\right\}\,,\quad
             l,l'\in\{0,\ms \frac12.,1,\ms \frac32.,\dots\}\,.
\]
This decomposition structure implies that
\[
 l\approx l'\quad\mr iff. \quad
         l,l'\;\; \mr are~both~integers~or~both~\mbox{half\,-\,integers}.\,.
\]
We can finally conclude that
\[
 \ot C.(\mr SU.\!(2))=\left\{[0],[\ms \frac12.]\right\}
                     \cong \Z_2\cong \wwh {\al C.(\mr SU.\!(}.2))\,.
\]
For other examples with finite and compact Lie groups see
\cite[Subsection~5.1]{Lledo04a}.
\end{eje}

We will now concentrate on the relation of the chain group $\ot C.(\al G.)$,
associated with the compact group $\al G.$ of a Hilbert C*-system $\HS$,
with the irreducible canonical endomorphisms restricted to $\al Z.$.
In particular recall the automorphisms on $\al Z.$ given in
Eq.~(\ref{alphaD}) by
$\alpha_\mt D.:=\rho_\mt D.\rest\al Z.\in\mr Aut.\al Z.$
which are associated with any class $D\in\wh{\al G.}$.
For a complete proof of the next theorem see 
Theorem~5.7 of \cite{Lledo04a}.

\begin{teo}\label{EndoChain}
\begin{itemize}
\item[(i)]
Let $D,D'\in\wh{\al G.}$ be equivalent, i.e.$~D\approx D'$. Then
$\alpha_\mt D.=\alpha_\mt D'.$ and we can
associate the automorphism $\alpha_\mt [D].\in\mr Aut.\al Z.$
with the chain group element $[D]\in\ot C.(\al G.)$.

\item[(ii)] There is a natural group homomorphism from the
chain group to the automorphism group
generated by the irreducible endomorphisms restricted to
$\al Z.$:
\begin{equation}\label{HomoCA}
\ot C.(\al G.)\ni [D]\mapsto \alpha_\mt [D].\in\mr Aut.\al Z.\,.
\end{equation}
\end{itemize}
\end{teo}

\begin{rem}\label{CompletePicture}
Note that the chain group and in particular Theorem~\ref{EndoChain}~(i)
completes the picture of the action of the irreducible canonical
endomorphisms on the center $\al Z.$ of the fixed-point algebra $\al A.$
(recall also Eq.~(\ref{alphaD})). Indeed, we may summarize this
action by means of the following diagram
\begin{eqnarray*}
  \wh{\al G.}&\to& \ot C.(\al G.)\;\; \to\;\; \mr Aut.\,\al Z. \\
    D        &\mapsto& \; [D] \;\;\;\mapsto\;\; \alpha_\mt [D].
\end{eqnarray*}
\end{rem}

\begin{teo}\label{GeneralZMap}
Let $\rho$ be a (reducible) canonical endomorphism.
Then its action on $\al Z.$ can be described
by means of the following formula
\[
 \rho(Z)=\sum_{[\mt D.]\in\ot C.(\al G.)} \alpha_\mt [D].(Z)\cdot
            E_{\mt {[D]}.}\;, \quad Z\in\al Z.\,,
\]
where $E_\mt D'.\in\al A.$ is the isotypical projection w.r.t.~$D'\in\wh{\al G.}$
and $E_{\mt {[D]}.}:=\mathop{\sum}\limits_{\mt {D'\in [D]}.} E_\mt D'.$.
For $n\geq 2$ we have
\[
 \rho^n(Z)=\sum_{[\mt {D_1}.],\dots,[\mt {D_n}.]}
              \alpha_\mt {[D_1\times\dots\times D_n]}.(Z)\cdot
              E_{\mt {[D_1]}.} \cdot\rho (E_{\mt {[D_2]}.})\cdot
              \ldots \cdot\rho^{n-1} (E_{\mt {[D_n]}.})
              \;, \quad Z\in\al Z.\,.
\]
\end{teo}
\begin{beweis}
The first equation is shown in \cite[Theorem~5.9]{Lledo04a}. From this
one can easily check the expression
for higher powers of $\rho\rest\al Z.$.
\end{beweis}

\section{Cuntz-Pimsner algebras}\label{Cuntz_Pimsner}
\label{sec_cp}

In the present section, we introduce some basic
properties of the Cuntz-Pimsner algebra
(also called Cuntz-Krieger-Pimsner algebra). The basic
reference for this topic is \cite{Pimsner97}. In the present paper
we will also use the alternative categorical
approach of \cite{Doplicher98}.

\subsection{Basic definitions}

Let $\ot R.$ be a {\it C*}-algebra, ${\ot H.}$ a countably generated Hilbert
$\ot R.$-bimodule. We assume that
${\ot H.}$ is full as a right Hilbert $\ot R.$-module.
We denote by $\al L. ( \ot H. , \ot H. )$ the {\it C*}-algebra
of adjointable, right $\ot R.$-module operators on ${\ot H.}$, and by
$\al K. ( \ot H. , \ot H. ) \subseteq \al L. ( \ot H. , \ot H. )$
the ideal of compact operators generated by the maps
\begin{equation}
\label{def_theta}
\theta_{\psi , \psi' } \in \al L. ( \ot H. , \ot H. )\;,\;\;
\psi , \psi'\in{\ot H.}\;,
\quad\mr with.\quad
\theta_{\psi , \psi' } (\varphi) :=
\psi \left \langle  \psi' , \varphi \right \rangle \ , \
\varphi \in {\ot H.} \ \ ,
\end{equation}
where $\left \langle \cdot , \cdot \right \rangle$ is the
$\ot R.$-valued scalar product defined on ${\ot H.}$. We also denote by
\[
\alpha\colon \ot R. \rightarrow \al L. ( \ot H. , \ot H. )
\]
the left $\ot R.$-action on ${\ot H.}$,
which we assume to be non-degenerate.

Let now $\ot R.$ be unital, and ${\ot H.}$ finitely generated as a right
Hilbert $\ot R.$-module by a subset $\Psi := \left\{ \psi_l \right\}_{l=1}^n$.
In this case, $\al L. ( \ot H. , \ot H. ) = \al K. ( \ot H. , \ot H. )$.
We consider the $*$-algebra $^0 {\al O.}_{\ot H.}$ generated by
a subalgebra *-isomorphic to $\ot R.$
and $\Psi$, subject to the relations{\footnote{With an abuse of notation, 
we identify $\ot R.$ with its
image in $^0 {\al O.}_{\ot H.}$.}}
\begin{eqnarray}
\psi^*_l \psi_m &=&\left \langle \psi_l , \psi_m \right\rangle
\label{def_pimsner}\\
A \psi_l     &=&  \alpha (A) \psi_l \;,
 \quad A \in \ot R.  \nonumber \\
 \sum_l \psi_l \psi_l^* &=& \1  \,. \label{eq_id}
\end{eqnarray}
Note that (\ref{def_pimsner}) implies
$\psi^* \psi' = \left \langle \psi , \psi' \right\rangle$ as well as
$\psi' \psi^* \varphi = \theta_{\psi' , \psi} (\varphi)$,
$\psi' , \psi^* , \varphi \in {\ot H.}$. Therefore one has the natural identification
\begin{equation}
\label{eq_compact_operators}
\theta_{\psi' , \psi} = \psi' \psi^* \ \ .
\end{equation}

The isomorphism class of $^0 {\al O.}_{\ot H.}$
does not depend on the choice of the set of generators
and it can be proven that there is a unique {\it C*}-norm on
$^0 {\al O.}_{\ot H.}$ such that the action of the circle
$\T:=\{z\in\C\mid |z|=1\}$ given by
\begin{equation}
\label{def_circle_action}
\delta\colon \T \rightarrow \aut ^0 {\al O.}_{\ot H.} \ \ , \ \
\delta_z (\psi) := z \psi \ \ , \ \
z \in \T , \ \psi \in {\ot H.} \ ,
\end{equation}
extends to an (isometric) automorphic action.
The Cuntz-Pimsner algebra ${\al O.}_{\ot H.}$ is by definition the
completion of $^0 {\al O.}_{\ot H.}$ w.r.t.~this norm. In
the case $\ot R. = \C$, ${\ot H.} = \C^d$, $d \in \N$, we
obtain the Cuntz algebra $\al O._d$. We denote by
\begin{equation}\label{SpecS1}
{\al O.}^k_{\ot H.} := \left\{
T \in {\al O.}_{\ot H.} \mid \delta_z (T) = z^k T \ , \ z \in \T
\right\}
\ , \
k \in \Z \ ,
\end{equation}
the spectral subspaces w.r.t. the circle action.

\begin{rem}\label{GenCP}
For the kind of Hilbert bimodules that we consider there is an
alternative way to generate the corresponding
Cuntz-Pimsner algebras that we will need in the following
sections. From \cite{PinzariIn97} we have that the spectral
subspaces are given as the following inductive limit, where
the natural inclusions are specified by tensoring from the
right by the identity on $\ot H.$:
\[
 ^0 {\al O.}_\ot H.^k = \lim_{\longrightarrow} \al L.(\ot H.^r,\ot H.^{r+k})\,.
\]
Moreover, we have the following natural $\Z$ grading
\[
 ^0 \al O._\ot H.={\oplus_{k\in\Z}}\;  ^0 \al O._\ot H.^k
\]
and ${\al O.}_\ot H.$ is the closure of $^0 \al O._\ot H.$
in the unique {\it C*}-norm for which $\delta$ is isometric.
\end{rem}

In order to discuss
universality properties of the Cuntz-Pimsner algebra, we give the
following definition:

\begin{defi} \label{def_hba} (\cite[\S 2]{Doplicher98})
Let $\ot R. \subset \ot B.$ be a {\it C*}-algebra inclusion. A {\bf Hilbert
$\ot R.$-bimodule in} $\ot B.$ is a closed vector space
${\ot H.} \subset \ot B.$,
such that $A \psi \in {\ot H.}$, $\psi A \in {\ot H.}$,
$\psi^* \psi' \in \ot R.$ for
every $A \in \ot R.$, $\psi , \psi' \in {\ot H.}$.
\end{defi}

We describe in explicit terms some well-known properties for the case of
Hilbert bimodule in {\it C*}-algebras. Our Hilbert bimodule $\ot H.$ is
finitely generated if there is a finite subset
$\left\{ \psi_l \right\}_{l=1}^n \subset {\ot H.}$
such that
$\psi = \sum_l \psi_l (\psi_l^* \psi)$, $\psi \in {\ot H.}$.
Moreover, $\ot H.$ is full if for every $A \in \ot R.$ there are
$\psi ,\psi' \in {\ot H.}$ such that $A = \psi^* \psi'$
and $\ot H.$ is non-degenerate if for any $A \in \ot R.$ with
$A \psi = 0$ for all $\psi \in {\ot H.}$, we have that $A = 0$.
If ${\ot H.}$ is finitely generated, then it is trivial to verify that
\[
P_{\ot H.} := \sum_l \psi_l \psi_l^* \ \in \ \ot B.
\]
is a projection, and that it does not depend on the choice of the set
of generators. We call $P_{\ot H.}$ the {\em support} of
${\ot H.}$ in $\ot B.$.

\begin{pro} \label{prop_uni} (\cite[Theorem~3.12]{Pimsner97})
Let $\ot R. \subset \ot B.$ be an inclusion of unital {\it C*}-algebras,
${\ot H.} \subset \ot B.$ a non-degenerate and full
Hilbert $\ot R.$-bimodule in $\ot B.$
with support $\1$. Then, there is a canonical monomorphism
${\al O.}_{\ot H.} \hookrightarrow \ot B.$.
\end{pro}

\begin{cor}
\label{cor_cp_hs}
Let $( \al F. , \al G. , \beta )$ be a Hilbert {\it C*}-system
with fixed-point algebra
$\al A. :=\al F.^\al G.$ and denote by $\al Z. :=\al A.' \cap \al A.$
its center.
Then, for every $D \in \widehat{\al G.}$ the Hilbert
$\al Z.$-bimodule $\ot H._\mt D. := \al Z. \al H._D$ induces a monomorphism
$\al O._{\ot H._D} \hookrightarrow \al F.$.
\end{cor}

\begin{rem} \label{rem_partial_isometry}
  \begin{itemize}
  \item[(i)]
  Examples of the above universality property can be found in the Cuntz-Pimsner
  algebra itself. Let $r = 1, 2 , \ldots $, and
  ${\ot H.}^r := {\ot H.} \otimes_{\ot R.} \cdots \otimes_{\ot R.} {\ot H.}$
  denote the $r$-fold tensor product with coefficients in $\ot R.$. Then, there is a
  natural identification
\[
{\ot H.}^r \simeq \left\{
\psi_1\cdot \ldots \cdot\psi_r \in {\al O.}_{\ot H.} \mid
\psi_k \in {\ot H.} , k = 1 , \ldots , r \right\} \ ,
\]
so that every ${\ot H.}^r$ is a Hilbert $\ot R.$-bimodule
 in ${\al O.}_{\ot H.}$, and
there are canonical morphisms $\al O._{{\ot H.}^r} \rightarrow {\al O.}_{\ot H.}$.
Note that if there is an element $S \in {\ot H.}^r$ such that
$\left \langle S , S \right \rangle = \1$, then $S$ appears
in the Cuntz-Pimsner algebra as an isometry, i.e.~$S^* S = \1$.
\item[(ii)]
The bimodule introduced in (\ref{GenBiMod}) is easily seen to be finitely generated,
full and non-degenerate.
\end{itemize}
\end{rem}

Let now $r,s \in \N$, and $\al K. ({\ot H.}^r , {\ot H.}^s)$ denote the set of compact,
right
$\ot R.$-module operators from ${\ot H.}^r$ into ${\ot H.}^s$. Applying
Eq.~(\ref{eq_compact_operators}), we may identify
\[
\al K. ({\ot H.}^r , {\ot H.}^s) \simeq
{\mathrm{span}}
\left\{
\psi_1 \cdots \psi_s\cdot \varphi_r^* \cdots \varphi_1^* \mid
\psi_h , \varphi_k \in {\ot H.} \,;\; h = 1 , \ldots , s \,;\; k = 1 ,\ldots , r
\right\} \ .
\]

To simplify notation we introduce the following
conventions: let ${\ot H.}$ be finitely generated by a set
$\Psi=\left\{ \psi_l \right\}_{l=1}^n$, $s \in \N$, and
$L := \left\{ l_1 , \ldots , l_s \right\} \in
\left\{ 1 , \ldots  ,n \right\}^s$
a multi-index of length $s$, i.e.~$|L| = s$. We denote
\begin{equation}
\label{def_tensor}
\psi_L := \psi_{l_1} \cdots \psi_{l_s} \in {\ot H.}^s \ ,
\end{equation}
so that by (\ref{eq_id}) we
find $$\sum_L \psi_L \psi_L^* = \1 \ ,$$ and we may also write
\begin{equation}
\label{eq_compact}
\al K. ({\ot H.}^r , {\ot H.}^s) \simeq
{\mathrm{span}}
\left\{
\psi_L \,A \,\psi_M^* \mid \;|L| = s \,;\; |M| = r \,;\; A \in \ot R.
\right\} \ .
\end{equation}

Finally, we need to introduce the notion of a nonsingular bimodule:
\begin{defi}
\label{def_nonsing}
A Hilbert $\ot R.$-bimodule ${\ot H.}$ with left action
$\alpha\colon \ot R. \rightarrow \al L. ( \ot H. , \ot H. )$
is called {\bf nonsingular} if
$\ \theta_{\psi , \psi} \in \alpha(\ot R.)$ for some $\psi \in {\ot H.}$
implies $\psi = 0$.
\end{defi}

\begin{pro}\label{pro_rel_com} (\cite[Proposition~3.5]{Doplicher98})
If ${\ot H.}$ is nonsingular and there is an isometry
$S \in {\ot H.}^n \cap \ot R.'$ for some $n \in \N$, $n > 1$,
then $C^*(S)' \cap {\al O.}_{\ot H.} = \ot R.$
(where $C^*(S)$
denotes the {\it C*}-algebra generated by $S$ in ${\al O.}_{\ot H.}$).
\end{pro}

\subsection{Endomorphisms in Cuntz-Pimsner algebras}
\label{EndCP}

It is well-known that a natural endomorphism is defined
over the Cuntz algebra
$\al O._n$, in the following way: if $\left\{ \psi_k \right\}_{k=1}^n$
is a set
of orthonormal isometries generating $\al O._n$, then we define
$\rho \in \endo \al O._n$ by
\begin{equation}
\label{def_can_end_1}
\rho (T) := \sum_k \psi_k T \psi_k^* \ \ , \ \ T \in \al O._n \ \ .
\end{equation}
Given a generic Hilbert $\ot R.$-bimodule ${\ot H.}$, it is not
possible to define
in a consistent way the analogue of (\ref{def_can_end_1}) over ${\al O.}_{\ot H.}$.
In fact,
the multiplicativity of $\rho$ is ensured by the fact that
$\psi_h^* \psi_k \in \al O._n \cap
\al O._n' = \C \1$, $h,k \in\{ 1 , \ldots , n\}$.
Given a set $\left\{ \psi_l \right\}_l$ of generators for
${\ot H.}$, it is not true in general that
$\psi_l^* \psi_m \in {\al O.}_{\ot H.} \cap {\al O.}_{\ot H.}'$.
Nevertheless, the
above obstruction for the existence of the endomorphism $\rho$
disappears if we consider a Hilbert bimodule which is free as a right
Hilbert $\ot R.$-module (cf.~Subsection~\ref{HCSsummary}), i.e.
\[
{\ot H.} \simeq \C^n \otimes \ot R. \,.
\]
In this case ${\ot H.}$ admits a (non-unique) set of generators
$\Psi := \left\{ \psi_k \right\}_{k=1}^n \subset {\ot H.}$ such that
$$\left \langle \psi_h , \psi_k \right \rangle = \delta_{hk}\1 \ .$$
We say in this case that $\Psi$ is a set of {\em
  orthonormal generators}.  Let $\al H._\Psi \subset \ot H.\subset
{\al O.}_{\ot H.}$ be the vector space spanned by elements of $\Psi$
(the latter regarded as a subset of ${\al O.}_{\ot H.}$). According to
the terminology of Definition~\ref{def_hba} we obtain by
(\ref{def_pimsner}) that $\al H._\Psi$ is a Hilbert $\C$-bimodule in
${\al O.}_{\ot H.}$. Such particular cases of Hilbert bimodules in
{\it C*}-algebras are called {\em Hilbert spaces in C*-algebras} in
\cite{Doplicher88}. As mentioned before we call them
{\em algebraic Hilbert spaces} (cf.~Section~\ref{HCSsummary} or
\cite{Lledo04a}).  It is clear that $\al H._\Psi$ depends on the
choice of $\Psi$. Different choices of orthonormal generating sets
$\Psi \subset {\ot H.}$ correspond to different isomorphisms ${\ot H.}
\simeq \al H._\Psi \otimes \ot R.$ of right Hilbert $\ot R.$-modules.
From the previous considerations, we obtain:

\begin{lem} \label{lem_sigmah}
  Let $\ot R.$ be a unital {\it C*}-algebra, ${\ot H.}$ a Hilbert $\ot
  R.$-bimodule, $\al H.$ a rank $n$ algebraic
  Hilbert space, $n \in \N$. Given a fixed
  isomorphism ${\ot H.} \simeq \al H. \otimes \ot R.$ of right
  Hilbert $\ot R.$-modules, there exists an endomorphism
  $\rho_\al H. \in \endo {\al O.}_{\ot H.}$, defined in Eq.~(\ref{def_can_end_1}),
  where $\left\{ \psi_k \right\}_{k=1}^n$ is a set of orthonormal
  generators corresponding to a basis of $\al H.$. The
  endomorphism $\rho_\al H.$
  does not depend on the choice of the basis of $\al H.$. Moreover, there is a
  unital monomorphism $\al O._n  \hookrightarrow {\al O.}_{\ot H.}$.
\end{lem}

\begin{rem}
Note that $\rho_\al H.$ commutes with the circle action, so that the spectral
subspaces introduced in (\ref{SpecS1}) are preserved,
i.e.~$\rho_\al H. ({\al O.}^k_{\ot H.} ) \subset {\al O.}^k_{\ot H.}$,
$k \in \Z$.
\end{rem}

\subsection{Amplimorphisms and their associated
             Cuntz-Pimsner algebras}
\label{sec_ampl}

Let $\ot R.$ be a unital {\it C*}-algebra with identity $\1$.
The {\it C*}-category
${\bf{bimod}} (\ot R.)$ having as objects the Hilbert
$\ot R.$-bimodules which are finitely generated and projective as right Hilbert
$\ot R.$-modules can be
described by means of {\em amplimorphisms},
i.e.~{\it C*}-algebra morphisms of the form
\[
\alpha \colon \ot R. \rightarrow \Mb{M}_n \otimes \ot R. \ \ ,
\]
\noindent where $\Mb{M}_n$ denotes the {\it C*}-algebra
of $n \times n$-matrices.
The Hilbert $\ot R.$-bimodule associated with $\alpha$ is defined by

\[
{\ot H.}_\alpha := \left\{
\psi \in \C^n \otimes \ot R. \mid \alpha (\1) \psi = \psi
\right\}
\]
with left and right $\ot R.$-actions defined as follows:
if $\psi := \left\{ A_l \right\}_{l=1}^n$,
$\psi' := \left\{ A'_m \right\}_{m=1}^n\in\ot H._\alpha$, we have
\[
\psi A := \left\{ A_l A  \right\}_l \in {\ot H.} \ \ , \ \
A \psi := \alpha (A) \psi \in {\ot H.} \ \ , \ \
\left \langle \psi , \psi' \right \rangle := \sum_l A^*_l A'_l \in \ot R. \ \ .
\]

The following notion of diagonal bimodule will play an important role
in the construction of examples presented in the following two sections.

\begin{defi}\label{AmplDiag}
A Hilbert $\ot R.$-bimodule ${\ot H.} = {\ot H.}_\alpha$,
$\alpha \in {\bf{ampl}}(\ot R.)$, is
called {\bf diagonal of rank $n$} if
$\alpha (A) := {\mathrm{diag}}(\alpha_1(A) , \ldots , \alpha_n(A))$,
$A \in \ot R.$, where $\alpha_1 , \ldots , \alpha_n \in \aut \ot R.$,
and ${\mathrm{diag}} (\cdot)$ denotes the diagonal matrix
in $\Mb{M}_n \otimes \ot R.$.
\end{defi}

It is clear that a diagonal bimodule ${\ot H.}$
is free as a right Hilbert $\al Z.$-module.
In this case the Cuntz-Pimsner algebra ${\al O.}_{\ot H.}$
is specified by the orthonormal basis
$\left\{ \psi_k \right\}_{k=1}^n \subset {\ot H.}$
that satisfies the following simple relations:
\begin{equation}
\label{eq_cp_diag}
\psi_h^* \psi_k = \delta_{hk} \1
\ \ , \ \
Z \psi_k = \psi_k \alpha_k (Z)
\ \ , \ \
\sum_k \psi_k \psi_k^* = \1
\ \ ,
\end{equation}
$h,k = 1 , \ldots , n$. Note that
 $\rho_{\al H.}(Z) \psi_k = \psi_k Z = \alpha^{-1}_k(Z) \psi_k$,
 $k = 1 , \ldots , n$. Let us consider the projections
$E_k := \psi_k \psi_k^* \in \al K. ( {\ot H.}  )$. Then
\begin{equation}
\label{eq_ce_z}
\rho_{\al H.} (Z) = \sum_k \alpha^{-1}_k(Z) E_k
\end{equation}

It will be useful for the construction of the class
of examples to give a characterization
of diagonal nonsingular bimodules
in the sense of Definition~\ref{def_nonsing} in order to apply
Proposition~\ref{pro_rel_com}. For this purpose we will apply
Gelfand's theorem and identify the unital abelian C*-algebra $\al Z.$
with continuous functions on the compact Hausdorff space $\Omega:=
\mathrm{spec}\,\al Z.$. Moreover, to any automorphism
$\alpha_k\in\mathrm{Aut}\,\al Z.$ there corresponds a homeomorphism $f_k$
of $\Omega$ such that for any $Z\in\al Z.$ we have the relation
\[
 \left(\alpha_k^{-1} \,Z\right)(\omega)=Z(f_k(\omega))\;,\quad \omega\in\Omega\;,
          k = 1 , \ldots , n\;.
\]

\begin{lem}
Let $\al Z.$ be an Abelian, unital {\it C*}-algebra and ${\ot H.} =
{\ot H.}_\alpha$ a diagonal free $\al Z.$-bimodule of rank $n$. If
$\left\{ \psi_k \right\}_{k=1}^n$ is an orthonormal  set of
generators for $\ot H.$ and $\varphi := \sum_h \psi_h Z_h \in {\ot
H.}$, then the equation $\theta_{\varphi , \varphi} = \alpha (Z)$
for some $Z\in\al Z.$ is equivalent to the equations
\begin{eqnarray}
Z_h Z_k^* &=& 0 \ \ , \ \ h \neq k \;,\quad h,k=1,\ldots, n\,,  \label{1nonsing}  \\
Z_k Z_k^* &=& \alpha_k (Z) \;,\quad k=1,\ldots, n\,. \label{2nonsing}
\end{eqnarray}
\end{lem}
\begin{beweis}
Evaluating the equation $\theta_{\varphi , \varphi} = \alpha (Z)$
on the basis elements we obtain
\[
 \sum_h \psi_h Z_h Z_k^* = \psi_k \alpha_k(Z)\;,\quad k=1,\ldots, n\,,
\]
which implies the statement.
\end{beweis}

\begin{pro}
\label{lem_nsing} Let $\al Z.$ be an Abelian, unital {\it
C*}-algebra and ${\ot H.} = {\ot H.}_\alpha$ a diagonal free $\al
Z.$-bimodule of rank $n$. Then, ${\ot H.}$ is nonsingular iff there
is a pair of indices $h,k \in \left\{ 1 , \ldots , n \right\}$, $h
\neq k$, such that $\alpha_h = \alpha_k \in \aut \al Z.$.
\end{pro}
\begin{beweis}
1.~Using Gelfand's theorem recall that to any automorphism $\alpha_i$ there
corresponds a homeomorphism $f_i$ on $\Omega:=\mathrm{spec}\,\al Z.$. Assume
that there exist a pair of indices
$h,k \in \left\{ 1 , \ldots , n \right\}$, $h \neq k$, such that
$f_h(\omega) = f_k(\omega)$, $\omega \in\Omega$. Suppose that there are
$Z,Z_1,\dots,Z_n\in\al Z.$ as above with $Z\not=0$ and such that
the equation $\theta_{\varphi , \varphi} = \alpha (Z)$ holds, i.e.
\begin{eqnarray}
Z_h(\omega) \overline{Z_k}(\omega) &=& 0 \ \ , \quad
                           \omega\in\Omega\;,\;h \neq k   \label{1nonsing-fn}  \\
|Z_k|^2\left(f_k(\omega) \right)&=& Z(\omega)
                \;,\quad h,k=1,\ldots, n\,. \label{2nonsing-fn}
\end{eqnarray}
Then for $\omega\in\Omega$ with $Z(\omega)\not=0$ we get from
Eq.~(\ref{2nonsing-fn}) that
\[
 Z_k(f_k(\omega))\not= 0 \not= Z_h(f_h(\omega))
\]
which contradicts Eq.~(\ref{1nonsing-fn}). Thus $Z=0$, hence $Z_k=0$, $k=1,\dots n$,
which implies that $\varphi=0$.
This shows that the bimodule is nonsingular.

2.~For the reverse implication we show that if there is a point
$\omega\in\Omega$ such that
for all pairs of indices $(h,k)$, $h\not=k$, we have $f_h(\omega)\not=f_k(\omega)$,
then $\ot H.$ is singular (i.e.~there exist $Z_1,\dots,Z_n\in\al Z.$ (not all equal
to 0) and a $Z\in\al Z.$ such that Eqs.~(\ref{1nonsing-fn}) and (\ref{2nonsing-fn})
hold; note that in this case it follows that $Z\not=0$).

We will show next that since $\Omega$ is a compact Hausdorff space
and the rank of the bimodule is finite
there is a neighborhood $W$ of $\omega$ such that
$f_k(W)\cap f_h(W)=\emptyset$ for all $(h,k)$, $h\not=k$. 
Indeed, for any pair $(h,k)$, $h \neq k$
since $f_h(\omega)\not=f_k(\omega)$ there exist open neighborhoods 
$W_{hk}$ resp.~$W_{kh}$ of $f_h(\omega)$ resp.~$f_k(\omega)$ with
$W_{hk}\cap W_{kh}=\emptyset$. Therefore we can define the 
open neighborhood of $\omega$ by
\[
 W:=\bigcap_{h\not= k} \left( f_h^{-1}(W_{hk}) \cap f_k^{-1}(W_{kh})\right)\;,
\]
which satisfies the required properties.

Let $\omega$ and $W$ be as in the preceding paragraph and define $Z(\cdot)$
as a positive continuous function with support contained in $W$.
Putting
\[
 Z_k(\omega):=\sqrt{Z\left(f_k^{-1}(\omega)\right)}
\]
we obtain continuous functions $Z,Z_1,\dots,Z_n\in\al Z.$ satisfying
Eqs.~(\ref{1nonsing-fn}) and (\ref{2nonsing-fn}). Hence $\ot H.$ is singular
and the proof is concluded.
\end{beweis}

\begin{rem}
\begin{itemize}
 \item[(i)] We note that if ${\ot H.}$ is a diagonal bimodule, then the left
       $\ot R.$-module action on ${\ot H.}$ is injective (i.e.~$\alpha (A) = 0
       \Rightarrow A = 0$, $A \in \ot R.$) and non-degenerate ($\alpha (\1) =
       \1_n$, where $\1_n$ is the identity of $\Mb{M}_n \otimes \ot R.$).
\item[(ii)] The Cuntz-Pimsner algebras $\al O._\ot H.$
       considered in this paper are nuclear. For a proof of this fact in the
       context of a more general class of Cuntz-Pimsner algebras see
       \cite{pLledo06}.

\end{itemize}
\end{rem}

\section{Examples of minimal {\it C*}-dynamical systems}\label{SecLie}

Let $\al G.$ be a compact group.
We denote by $\widehat{\al G.}$ its dual object
and by $\ot C.(\al G.)$ its chain group (cf.~Subsection~\ref{SecChain}).
Given a unital (not necessarily Abelian)
{\it C*}-algebra $\ot R.$, we consider a fixed chain
group action
\[
\ot C.(\al G.)\ni [D]\mapsto \alpha_\mt [D].\in\mr Aut.\ot R.\,.
\]
Later on we will restrict to the case
where $\ot R.=\al Z.$ is Abelian.

Let $V$ be a unitary representation of $\al G.$ on
a finite-dimensional Hilbert space $\al H.$.
In general, $V$ may be reducible, hence we can decompose
it as a direct sum
\[
 V= \sum_{\mt {D_i}.} E_{\mt {D_i}.}\, V\cong \oplus_i m_i\, V_{\mt {D_i}.}
  \quad \mr on.\quad
 \al H. = \sum_{\mt {D_i}.} E_{\mt {D_i}.}\, \al H.
          \cong \oplus_i m_i\, \al H._{\mt {D_i}.}\,,
\]
where
$m_i$ is the multiplicity of the irreducible representation
$V_{\mt {D_i}.}$ in $V$ and $E_\mt {D_i}.\in (V,V)$ is the
isotypical projection corresponding to $D_i\in\widehat{\al G.}$.
A useful orthonormal basis adapted to the previous decomposition
is given by
\begin{equation}\label{AdaptedONB}
 \Psi:=\left\{ \psi_{\mt {D,l,i}.}\mid D\in\widehat{\al G.}\,;\;
        i=1,\dots, d=\mr dim. D\,;\; l=1,\dots, m_\mt D.
 \right\}\,.
\end{equation}

Next we consider the free right $\ot R.$-module $\ot H.$
generated $\ot R.$ and $\Psi$ as in Eq.~(\ref{GenBiMod}).
Moreover, $\ot H.$ becomes a bimodule if we define the left action
of $\ot R.$ as
\begin{equation}\label{left_action}
 B\,\psi_{\mt {D,l,i}.}:=\psi_{\mt {D,l,i}.}\,\alpha_{\mt {[D]}.}^{-1}(B)\,,
\end{equation}
where $[D]\in\ot C.(\al G.)$ is the chain equivalence class
corresponding to $D\in\widehat{\al G.}$ (cf.~Subsection~\ref{SecChain}).
By construction $\ot H.$ is a diagonal bimodule in the sense
of Definition~\ref{AmplDiag}.

Let $\al O._\ot H.$ be the Cuntz-Pimsner algebra generated by
$\Psi$ and $\ot R.$ (cf.~Section~\ref{Cuntz_Pimsner}). Then,
$\al H.$ is an algebraic Hilbert space in $\al O._\ot H.$
with support $\1$. With respect to the basis (\ref{AdaptedONB}),
the isotypical projections may be written as
\begin{equation}\label{ED}
 E_\mt D.=\sum_{l,i}\psi_{\mt {D,l,i}.}\psi_{\mt {D,l,i}.}^*
\end{equation}
and since the support of $\al H.$ is $\1$ we also have the relation
\[
 \sum_{D,l,i}\psi_{\mt {D,l,i}.}\psi_{\mt {D,l,i}.}^*=\1\,.
\]
From Eq.~(\ref{ED}), we have that
$E_\mt D. B=B E_\mt D.$ as well as $E_\mt [D]. B=B E_\mt [D].$
(recall from Theorem~\ref{GeneralZMap} that
 $E_{\mt {[D]}.}:=\mathop{\sum}\limits_{\mt {D'\in [D]}.} E_\mt D'.$).

\begin{rem}
If the representation $V$ is irreducible
(e.g.,~in the case of the defining representation of SU($N$)), then for any
$T \in \al K. (\ot H.)$ we have simply
\begin{equation}\label{irredleftout}
T (B\psi) = B\,T(\psi)\,, \quad B \in \ot R.\,, \psi \in \ot H.\,.
\end{equation}
In fact, if $V\in D\in\widehat{\al G.}$, then we have for the
generators
\[
 T(B\psi_i)=T(\psi_i \alpha_\mt [D].^{-1}(B))
            =T(\psi_i) \alpha_\mt [D].^{-1}(B)
            = B T(\psi_i)\,.
\]
Note that Eq.~(\ref{irredleftout}) is no longer true
in the reducible case.
\end{rem}

To define an action of $\al G.$ on $\al O._\ot H.$ it is enough
to specify it by means of the representation $V$ on the
generating module $\ot H.$:
\begin{equation}
\label{eq_gact}
 g(\psi B):= (V(g)\psi) B
 \;,\quad \psi \in \al H.\,,\;B\in\ot R. \,,\; g \in \al G.\,.
\end{equation}
Note that this immediately implies that $g(B\psi)=  Bg(\psi)$,
$\psi \in \al H.$, $B\in\al B.$, in fact
\begin{equation}
\label{eq_bg}
g  (B \psi_{\mt {D,l,i}.}) =
g ( \psi_{\mt {D,l,i}.} \alpha_{\mt [D].}^{-1} (B) ) =
( V(g) \psi_{\mt {D,l,i}.} ) \alpha_{\mt [D].}^{-1} (B) ) =
B g ( \psi_{\mt {D,l,i}.} )
\ .
\end{equation}
Then $g$ extends to an automorphism $g \in \mr Aut.\al O._\ot H.$ and
$(\al O._\ot H.,\al G.)$ becomes a {\it C*}-dynamical
system. We denote its fixed-point algebra by
\[
 \al A.:=(\al O._\ot H.)^\al G.=\{A\in\al O._\ot H.\mid
           g(A)=A\,,\; g\in\al G.\}\,.
\]

Let now $\rho_\al H. \in \endo {\al O.}_{\ot H.}$ be the
endomorphism induced by an orthonormal basis
$\left\{ \psi_i \right\}_i^d \subset \al H.$,
according to Lemma~\ref{lem_sigmah}.
As we mentioned in Subsection~\ref{HCSsummary}
we call canonical endomorphism the restriction of
$\rho_\al H.$ to $\al A.$, i.e.
\[
\rho := \rho_{\al H.} \rest {\al A.} \in \endo \al A. \ .
\]

\begin{lem} \label{lem_AZrho}
 Let $\al A.$ be the fixed-point algebra of
 the {\it C*}-dynamical system $(\al O._\ot H.,\al G.)$ and
 $\rho$ be the endomorphism introduced before. Then
  \begin{itemize}
  \item[(i)] $\rho^n( \ot R.)
             :=\left\{ \rho^n(B) \mid B \in \ot R.\right\}\subseteq\al A.$,
             $n\in\N$.
  \item[(ii)] $\rho(B) B'= B' \rho(B)$, $B,B'\in\ot R.$.
  \end{itemize}
\end{lem}
\begin{beweis}
By the definition of the action of $\al G.$ on the Cuntz-Pimsner algebra
specified in Eq.~(\ref{eq_gact}), it follows that $\ot R.\subseteq\al A.$.
Since $\rho$ is an endomorphism of $\al A.$ it follows that
$\rho^n( \ot R.)\subseteq\al A.$.
The proof of part (ii) follows immediately from the definition of 
$\rho$ and Eq.~(\ref{left_action}).
\end{beweis}

\begin{lem}\label{Decomp1}
Let $V$ be a representation of the compact group
$\al G.$ acting on the algebraic Hilbert space $\al H.$. Then,
the trivial representation $\iota$ is contained
in the decomposition of $\mathop{\otimes}\limits^n V$,
i.e.~$\iota\prec\mathop{\otimes}\limits^n V$,
iff there exists an isometry $S \in \al H.^n$ such that
$g(S)=S$, $g\in\al G.$, (i.e.~$S \in \al H.^n \cap \al A.$ ).
\end{lem}
\begin{beweis}
The unitary representation $\mathop{\otimes}\limits^n V$ contains
the trivial representation if and only if there exists a
(normalized) vector $S \in \al H.^n$, invariant under
$\mathop{\otimes}\limits^n V$.
Now, by Remark~\ref{rem_partial_isometry}~(i), $S$ appears as an
isometry in ${\al O.}_{\ot H.}$. Moreover, by the
definition of the $\al G.$-action we have $g(S)=S$,
thus $S \in \al A. \cap \al H.^n$. Conversely, if $S \in \al H.^n$
is a $\al G.$-invariant isometry, then $\al H._\iota:=\C S\subset\al H.^n$
carries the trivial representation, hence
$\iota\prec\mathop{\otimes}\limits^n V$.
\end{beweis}

\begin{lem}\label{XnRho}
Let $X \in \al H.^n$. Then, the equation $X B = \rho^n(B) X$ holds for all
$B\in\ot R.$.
\end{lem}
\begin{beweis}
Consider $X = \varphi_1\cdot\ldots\cdot \varphi_n\in\al H.^n$,
$\varphi_k \in \al H.$, $k=1,\dots,n$.
Each $\varphi_k$ can be decomposed in terms of the chain group components,
i.e. $\varphi_k=\sum_{\mt {[D_{i_k}]}.} \varphi_{\mt {[D_{i_k}]}.}$,
where $\varphi_\mt {[D_{i_k}]}. := E_\mt {[D_{i_k}]}. \varphi_k$ and
$E_{\mt {[D]}.}:=\mathop{\sum}\limits_{\mt {D'\in [D]}.} E_\mt D'.$
(cf.~Theorem~\ref{GeneralZMap}).
Therefore we have
\[
 X=\sum_{\mt {[D_{i_1}],\dots,[D_{i_n}]}.}
   \varphi_{\mt {[D_{i_1}]}.}
   \cdot\ldots\cdot \varphi_{\mt {[D_{i_n}]}.}\,.
\]

From the definition of the left action
in Eq.~(\ref{left_action}) we have for any $B\in\ot R.$
\begin{eqnarray*}
XB &=& \sum_{\mt {[D_{i_1}],\dots,[D_{i_n}]}.}
       \varphi_{\mt {[D_{i_1}]}.}\cdot\ldots\cdot \varphi_{\mt {[D_{i_n}]}.}\cdot B \\
   &=& \sum_{\mt {[D_{i_1}],\dots,[D_{i_n}]}.}
       \alpha_\mt {[D_{i_1}\times\dots\times D_{i_n}]}.(B)
       \varphi_{\mt {[D_{i_1}]}.}\cdot\ldots\cdot \varphi_{\mt {[D_{i_n}]}.} \\
   &=& \Big(\alpha_\mt {[D_{i_1}\times\dots\times D_{i_n}]}.(B)
       \cdot E_{\mt {[D_1]}.} \cdot\rho (E_{\mt {[D_2]}.})\cdot
       \ldots \cdot\rho^{n-1} (E_{\mt {[D_n]}.}) \Big)
       \varphi_1\cdot\ldots\cdot \varphi_n       \\
   &=& \rho^n(B) \,X\,,
\end{eqnarray*}
where for the last equation we have used Theorem~\ref{GeneralZMap}.
\end{beweis}

\begin{pro}
\label{pro_rc}
  Let $X\in\al H.^n$. If $X$ is $\al G.$-invariant,
  i.e.~$X\in\al A.\cap\al H.^n$, then
  $XB=BX$, $B\in\ot R.$.
\end{pro}
\begin{beweis}
Let $X=\varphi_1\cdot\ldots\cdot \varphi_n\in\al H.^n$,
$\varphi_i\in\al H.$, $i=1,\dots,n$.
First note that, as in the proof of Lemma~\ref{XnRho}, we have
\begin{equation}\label{PasaZ}
XB= \sum_{\mt {[D_{i_1}],\dots,[D_{i_n}]}.}
    \alpha_\mt {[D_{i_1}\times\dots\times D_{i_n}]}.(B)
    \varphi_{\mt {[D_{i_1}]}.}\cdot\ldots\cdot \varphi_{\mt {[D_{i_n}]}.}
     \;,\quad B\in\ot R..
\end{equation}
Moreover, from the definition of the isotypical projections in
Eq.~(\ref{ED}) it is clear that $g(E_{\mt {D_k}.})=E_{\mt {D_k}.}$,
$k=1,\dots,n$. Therefore, denoting $\varphi_{\mt {D_{i_k}}.}
:=E_{\mt {D_{i_k}}.}\varphi_k$ we have that
\begin{equation}\label{ChainCompInv}
 g(X) = X\;,\;g\in\al G., \quad \mathrm{iff} \quad
        g\left(\varphi_{\mt {D_{i_1}}.}
        \cdot\ldots\cdot \varphi_{\mt {D_{i_n}}.}\right)
      = \varphi_{\mt {D_{i_1}}.}\cdot\ldots\cdot
        \varphi_{\mt {D_{i_n}}.}\;,\;g\in\al G.\,,
\end{equation}
for all $D_{i_k}$ appearing in the decomposition of $V$.
(Note that the preceding equation implies in particular the invariance of the
corresponding chain group components, i.e.~
$g(\varphi_{\mt {[D_{i_1}]}.}
        \cdot\ldots\cdot \varphi_{\mt {[D_{i_n}]}.})
      = \varphi_{\mt {[D_{i_1}]}.}\cdot\ldots\cdot
        \varphi_{\mt {[D_{i_n}]}.}$.)
Now, decomposing the tensor product
$\al H._{\mt {D_{i_1}}.}\cdot\ldots\cdot\al H._{\mt {D_{i_n}}.}$
in terms of irreducible components it is clear that Eq.~(\ref{ChainCompInv})
implies that there is a $\al G.$-invariant isometry
carrying the trivial representation,
i.e.~$\iota\in D_{i_1}\times\dots\times D_{i_n}$
(cf.~Lemma~\ref{Decomp1}). By the definition
of the chain group (cf.~Definition~\ref{DefChain}) this shows that
\[
\alpha_\mt {[D_{i_1}\times\dots\times D_{i_n}]}.(B)
 =\alpha_\mt [\iota].(B)=B.
\]
Therefore Eq.~(\ref{PasaZ}) reads $XB=BX$, $B\in\ot R.$.
\end{beweis}

For the rest of the present section, we assume that
the coefficient algebra $\ot R.$ is Abelian
and will write $\al Z.$ instead of $\ot R.$.

\begin{pro}
\label{teo_hrs}
Let $T \in \al L. ( \ot H.^s , \ot H.^n)$, $s,n\in\N$.
If $T$ is $\al G.$-invariant, i.e.~
$T \in \al A.\cap\al L. ( \ot H.^s , \ot H.^n)$,
then $TZ=ZT$, $Z \in \al Z.$. Moreover, $\al Z.\subseteq \al A.'\cap\al A.$.
\end{pro}
\begin{beweis}
We divide the proof in two steps. First,
we consider an elementary tensor of the form
\[
X = \varphi_1 \cdot \ldots \cdot \varphi_n \cdot
    \tilde{\varphi}_1^* \cdot \ldots \cdot \tilde{\varphi}_s^*
    \;\in \;\al H.^n (\al H.^*)^s \;\subset\; \al L. ( \ot H.^s , \ot H.^n )
\]
where $\varphi_i  , \tilde{\varphi_k} \in \al H.$, $i, = 1 \ldots n$,
$k = 1 , \ldots , s$. Recall that
$\al H.^n (\al H.^*)^s$ carries a unitary
representation of $\al G.$ equivalent to
$V^{\otimes^n} \otimes {\overline V}^{\otimes^s}$.
Decomposing again each $\varphi_i$ in terms of the chain group components
as in the proof of Lemma~\ref{XnRho} we get
\[
X=\sum_{\substack{ {\mt {[D_{i_1}],\dots,[D_{i_n}]}.}\\
       {\mt {[\overline{D}_{k_1}],\dots,[\overline{D}_{k_s}]}.}}}
       \varphi_{\mt {[D_{i_1}]}.}\cdot\ldots\cdot
       \varphi_{\mt {[D_{i_n}]}.}\cdot
       \varphi_{\mt {[\overline{D}_{k_1}]}.}\cdot\ldots\cdot
       \varphi_{\mt {[\overline{D}_{k_s}]}.}\,,
\]
where $\varphi_\mt {[D_{i_k}]}. := E_\mt {[D_{i_k}]}. \varphi_k$.
Therefore
\[
 XZ= \sum_{\substack{ {\mt {[D_{i_1}],\dots,[D_{i_n}]}.}\\
          {\mt {[\overline{D}_{k_1}],\dots,[\overline{D}_{k_s}]}.}}}
    \alpha_\mt {[D_{i_1}\times\dots\times D_{i_n}\times
                \overline{D}_{k_1}\times\dots\times
                 \overline{D}_{k_s}]}.(Z)
    \varphi_{\mt {[D_{i_1}]}.}\cdot\ldots\cdot
    \varphi_{\mt {[\overline{D}_{k_s}]}.}\,.
\]

Suppose now that $g(X)=X$, $g \in \al G.$. As in the proof of
Proposition~\ref{pro_rc} this implies that the product of chain
group components are also $\al G.$-invariant (cf.~Eq.~(\ref{ChainCompInv})).
Therefore, Lemma~\ref{Decomp1} implies that
$\iota\in D_{i_1}\times\dots\times D_{i_n}\times
 \overline{D}_{k_1}\times\dots\times\overline{D}_{k_s}$, hence
\[
\alpha_\mt {[D_{i_1}\times\dots\times D_{i_n}\times
                \overline{D}_{k_1}\times\dots\times
                 \overline{D}_{k_s}]}.(Z)
= Z \,.
\]
This shows that $XZ = ZX$, $Z\in\al Z.$.

Second, we consider an element of the full intertwiner space
$T \in \al L. ( \ot H.^s , \ot H.^n)$.
Using the identification in (\ref{eq_compact}), we can write
$T$ in terms of an orthonormal basis $\left\{ \psi_i \right\}
\subset\al H.$ as follows
\[
T = \sum_{I,J} \psi_I\, T_{IJ}\, \psi_J^* \ ,
\]
where, to simplify notation, we put
$\psi_I := \psi_{i_1} \cdot\ldots\cdot \psi_{i_n} \in \al H.^n$,
$\psi_J := \psi_{j_1} \cdot\ldots\cdot \psi_{j_s} \in \al H.^s$, and
\[
T_{IJ} = \psi_I^* \,T\, \psi_J \in \al Z. \ \ .
\]
Thus, by Lemmas~\ref{XnRho} and \ref{lem_AZrho}~(i) we have
\[
T = \sum_{I,J} \psi_I\, \psi_J^* \,\rho^s (T_{IJ}) \,
 \quad\mr with. \quad \rho^s (T_{IJ})\in\al Z. \,.
\]
Suppose now that $g (T) = T$, $g \in \al G.$.
Using the group mean given by
\begin{equation}\label{group-mean}
\ot m._\al G.(H):=\int_{\al G.} g(H)\ dg\;,\quad H\in\al O._\ot H.\,,
\end{equation}
where $dg$ is the normalized Haar measure, we obtain
\[
T = \ot m._\al G.(T)
  =\sum_{I,J} X_{IJ} \rho^s ( T_{IJ} ) \,,
\]
where $X_{IJ} := \int_{\al G.}  g (\psi_{I} \psi_J^*) \, dg \,
\in\al H.^n (\al H.^*)^s$. By definition the $X_{IJ}$
are $\al G.$-invariant, hence the first part of the proof
gives $X_{IJ}Z = ZX_{IJ}$, $Z \in \al Z.$.
Applying Lemma~\ref{lem_AZrho}~(ii) we obtain finally
\[
ZT = \sum_{I,J} ZX_{IJ} \rho^s ( T_{IJ} ) =
\sum_{I,J} X_{IJ} \rho^s ( T_{IJ} ) Z =
TZ \;,\quad Z\in\al Z.\,.
\]

To show the inclusion $\al Z.\subseteq\al A.'\cap\al A.$ recall
that from the definition of the group action in Eq.~(\ref{eq_gact})
we already have $\al Z. \subseteq \al A.$. Next we show that
$\left\{  \al L. ( \ot H.^s , \ot H.^n) \cap \al A.  \right\}_{s,n \in \N}$
is dense in $\al A.$: first note that since the circle action given
in Eq.~(\ref{def_circle_action}) commutes with the group action, we
may decompose $\al A.$ in the corresponding spectral subspaces
(cf.~(\ref{SpecS1})),
\[
 \al A.=\oplus_{k\in\Z}\al A._k\,.
\]
The density follows from the existence of a sequence of
norm one projections $E_r$
from $\al A._k$ onto $\al L.(\ot H.^r,\ot H.^{r+k})\cap\al A.$
pointwise convergent to the identity
(cf.~the proof of Proposition~3.4~(b) in \cite{Doplicher98}).
Therefore, by the first part of this theorem any $Z\in\al Z.$ will commute
with any $T\in\al L. ( \ot H.^s , \ot H.^n) \cap \al A.$, $s,n\in\N$,
hence with any $A\in\al A.$. This shows that
$\al Z.\subseteq\al A.'\cap\al A.$.
\end{beweis}

In this section we have shown that given a chain group action
$\alpha  \colon \ \ot C.(\al G.) \to \mr Aut.(\al Z.)$ we
can construct for any finite-dimensional representation $V$ of
the compact group $\al G.$ a diagonal Hilbert bimodule $\ot H.$ and
a {\it C*}-dynamical system $(\al O._\ot H.,\al G.)$ with the
corresponding Cuntz-Pimsner algebra. To show the minimality of
$(\al O._\ot H.,\al G.)$, i.e.~$\al A.'\cap\al O._\ot H.=\al Z.$,
where $\al A.=\al O._\ot H.^\al G.$ is the corresponding fixed point
algebra, we need to make some further assumptions on the representation
$V$. Essential facts in the proof of the minimality property are
that $\ot H.$ is nonsingular and that there exists a $\al G.$-invariant
isometry (recall Definition~\ref{def_nonsing} and
Proposition~\ref{pro_rel_com}). These facts will be guaranteed if we
consider finite-dimensional representations in the following class
(see also Proposition~\ref{lem_nsing}):

\begin{defi}
\label{def_gfin}
Let $\al G.$ be a compact group. We denote by $\gfin$
the set of all finite-dimensional representations $V$ of $\al G.$ satisfying the
following properties:
\begin{itemize}
\item[(i)] $V$ has an irreducible subrepresentation of dimension
           or multiplicity $\geq 2$.
\item[(ii)] There exists an $n \in \N$ such that
            $\iota \prec \mathop{\otimes}\limits^n V$.
\end{itemize}
\end{defi}

\begin{rem}
If $\al G.$ is non-Abelian we can use as representation $V$ any irreducible
representation of dimension $\geq 2$. If $\al G.$ is Abelian we consider a
representation containing some character with multiplicity $\geq 2$.

Part (ii) in the preceding definition is satisfied if the representation
$V$ on the Hilbert space $\al H.$ satisfies det$V=1$. In this case
there exists an isometry $S \in \al H.^d$, $d=\mr dim.\al H.$, with
$g (S) = S$, $g \in \al G.$. We may pick $S$ as a normalized vector
generating the totally antisymmetric tensor power
$\wedge^d \al H. \subseteq \al H.^d$. If $V$ does not have determinant $1$
we can always consider $V\oplus \overline{\mr det.V}$.
\end{rem}

\begin{teo}
\label{thm_main}
Let $\al G.$ be a compact
group, $\al Z.$ a unital Abelian C*-algebra
and consider a fixed chain group action
$\alpha  \colon \ \ot C.(\al G.) \to \mr Aut.(\al Z.)$.
For any $V\in\gfin$ there exists a nonsingular diagonal bimodule
$\ot H._\mt V.$ with left $\al Z.$-action given in terms
of $\alpha$ as in Eq.~(\ref{left_action}) and
a C*-dynamical system $(\al O._{\ot H._\mt V.} ,\al G.)$
satisfying the following properties:
\begin{itemize}
\item[(i)] $(\al O._{\ot H._\mt V.} ,\al G.)$ is minimal,
i.e.~$\al A._\mt V.' \cap \al O._{\ot H._\mt V.} = \al Z.$, where
$\al A._\mt V. := \al O._{\ot H._\mt V.}^{\al G.}$
is the corresponding fixed-point algebra.
\item[(ii)] The Abelian C*-algebra $\al Z.$ coincides with
the center of the fixed-point algebra $\al A._\mt V.$, i.e.~
$\al A._\mt V.' \cap \al A._\mt V. = \al Z.$.
\end{itemize}
Moreover, if $\al G.$ is a compact Lie group, then the Hilbert
spectrum of $(\al O._{\ot H._\mt V.}  ,\al G.)$ is full, i.e.~for each
$D\in\wh{\al G.}$ there is an invariant algebraic Hilbert space
$\al H._{\mt D.} \subset\al O._{\ot H._\mt V.}$
(in this case not necessarily of support $\1$)
such that $\al G.\restriction \al H._{\mt D.} \in D$.
\end{teo}
\begin{beweis}
Consider a representation $V\in\gfin$ and let
$\alpha  \colon \ \ot C.(\al G.) \to \mr Aut.(\al Z.)$
be a fixed chain group action. We construct
the diagonal Hilbert $\al Z.$-bimodule $\ot H._\mt V.$
with a left action as specified in Eq.~(\ref{left_action}).

First we show that $\ot H._\mt V.$ is nonsingular.
By part (i) in Definition~\ref{def_gfin}
there are indices
$i\not=j$ such that for the corresponding orthonormal basis elements
adapted to the decomposition of $V$ (cf.~Eq.~(\ref{AdaptedONB}))
satisfy $\psi_{\mt {D,l,i}.}\not=\psi_{\mt {D,l,j}.}$
as well as
\[
 Z\psi_{\mt {D,l,i}.}=\psi_{\mt {D,l,i}.}\alpha_{\mt [D].}^{-1}(Z)
 \quad\mr and.\quad
 Z\psi_{\mt {D,l,j}.}=\psi_{\mt {D,l,j}.}\alpha_{\mt [D].}^{-1}(Z)\;,
 \;\; Z\in\al Z.\,.
\]
From Proposition~\ref{lem_nsing} we conclude that $\ot H._\mt V.$
is nonsingular.

We define the {\it C*}-dynamical system
$( {\al O.}_{\ot H._\mt V.}  , \al G. )$ as in (\ref{eq_gact}).
Next we show that $( {\al O.}_{\ot H._\mt V.}  , \al G. )$ is minimal
(cf.~Definition~\ref{defs2-1}). By part (ii) in Definition~\ref{def_gfin}
and Proposition~\ref{pro_rel_com}
we conclude that $C^*(S)' \cap {\al O.}_{\ot H._\mt V.} = \al Z.$.
But $C^*(S) \subset \al A._\mt V.$, thus
$\al A._\mt V.' \cap {\al O.}_{\ot H._\mt V.} \subseteq \al Z.$.
The reverse inclusion follows from
$\al Z.\subseteq\al A._\mt V.'\cap\al A._\mt V.$ 
(cf.~Proposition~\ref{teo_hrs}).
This proves that $\al A._\mt V.' \cap {\al O.}_{\ot H._\mt V.} = \al Z.$,
i.e.~$( {\al O.}_{\ot H._\mt V.}  , \al G. )$ is minimal.
Note that the preceding equation also shows that
$\al A._\mt V.' \cap \al A._\mt V. = \al Z.$, since
$\al A._\mt V.' \cap \al A._\mt V. \subseteq
\al A._\mt V.' \cap {\al O.}_{\ot H._\mt V.} =\al Z.$
(the reverse inclusion $\al Z. \subseteq \al A._\mt V.' \cap \al A._\mt V.$ is
proved in Proposition~\ref{teo_hrs}).

If $\al G.$ is a compact Lie group, then there is a
faithful representation $V$ such that every irreducible
is contained in some tensor power $\otimes^n V$, $n\in\N$
(see \cite[Theorem~III.4.4]{bBroecker85}).
Therefore $(\al O._{\ot H._\mt V.},\al G.)$ has full Hilbert spectrum.
\end{beweis}

\begin{eje}\label{SUN} {\bf (Minimal {\it C*}-dynamical systems for SU$(N)$)}\\
Theorem~\ref{thm_main} can be applied to the group $\al G.=\mr SU.\!(N)$.
One can use the defining representation $V$ of dimension $N\geq 2$, since
it is easy to see that $V\in\gfin$. In fact, in this case we have
\[
 \iota \prec \mathop{\otimes}\limits^\mt N. V
\]
(see also Example~\ref{SU2_v1}). Therefore the preceding theorem gives a
minimal {\it C*}-dynamical system
$({\al O.}_{\ot H._\mt V.},\mr SU.\!(N))$ with
full spectrum.
\end{eje}

\section{Construction of Hilbert {\it C*}-systems}
\label{SecCompact}

In the present section we will construct Hilbert {\it C*}-systems
$(\al F.,\al G.)$ for compact (non-Abelian) groups $\al G.$.
This means that the {\it C*}-dynamical system must
encode the categorical structure of the dual of $\al G.$ (cf.~Definition~\ref{defs2-1}).
The important new feature now is that in $\al F.$ all irreducible
representations must be realized on algebraic Hilbert spaces with
support $\1$. As a first step towards this goal
we will construct a {\it C*}-algebra $\al F.$ suitably
generated by the Cuntz-Pimsner algebras $\al O._{\ot H._V}$, where $V\in\ot G.$
and $\ot G.$ is a family of finite-dimensional unitary representations
of $\al G.$.

\subsection{The {\it C*}-algebra of a chain group action}
\label{5.1}

Let $\al G.$ be a compact group with chain group $\ot C.(\al G.)$.
We consider a (non necessarily Abelian) {\it C*}-algebra $\ot R.$ with
identity $\1$ and carrying an action of the chain group
\[
\alpha\colon  \ot C.(\al G.) \rightarrow \aut \ot R. \ .
\]
With these ingredients, and a family
$\ot G.$ of finite-dimensional unitary
representations of $\al G.$,
we will construct a {\it C*}-dynamical system generalizing both Cuntz-Pimsner
algebras of the type given in
Section~\ref{sec_ampl} and crossed-products by Abelian groups
(cf.~\cite{Pimsner97}).

The following variation in the notation of the basis elements in this section
is convenient for the presentation of the algebra $\al F.$ given below.
Let $V \in \ot G.$ and consider the decomposition of 
$V$ into irreducible components. Then, we construct an orthonormal basis
for $\al H._\mt V.$, denoted by
\[
\Psi_\mt V. 
:= 
\left\{  
 \psi_{\mt V.,i} \ , \  i = 1 , \ldots , \mathrm{dim}\al H._V  
\right\}
\ ,
\]
and such that each $\psi_{\mt V.,i}$ transforms according an irreducible
subrepresentation of $V$ with class $D$ (see also Section~\ref{SecLie}):
\[
V(g) \psi_{\mt V.,i}
=
V_{\mt D.}(g)\psi_{\mt V.,i}
\ \ , \ \
g \in \al G.
\ .
\]
When $\psi_{\mt V.,i}$ satisfies the preceding condition, we write
\[
\psi_{\mt V.,i;(\mt D.)} := \psi_{\mt V.,i}
\ .
\]
We emphasize that in this case we are counting the basis elements in 
a different way as in Eq.~(\ref{AdaptedONB}). In fact, now
the element $D$ of the dual, which is written in brackets, 
does not play the role of an index.  
This label is used simply to recall the transformation character
of the basis element under
the action of $\al G.$.

We define a *-algebra $\ ^0 \al F. =$ $\ ^0 \al F. (\alpha , \ot G.)$,
generated by $\ot R.$ and  $\left\{ \Psi_\mt V. \right\}_{\mt {V\in \ot G.}.}$,
and satisfying the relations
\begin{eqnarray}
\sum_i \psi_{\mt V.,i;(\mt D.)} \ \psi_{\mt V.,i;(\mt D.)}^*
               &=&  \1   \label{31.1}    \\
\psi_{\mt V.,i;(\mt D.)}^* \ \psi_{\mt V.,j;(\mt D.')}
               &=& \delta_{ij} \1  \label{31.2} \\
\psi_{\mt V.,i;(\mt D.)} B
               &=&  \alpha_{\mt [D].} (B) \ \psi_{\mt V.,i;(\mt D.)}
                   \;,\; [D]\in\ot C.(\al G.)
                   \;,\; B\in\ot R.  \label{31.3}   \\
\left[ \psi_{\mt V.,i;(\mt D.)} \ , \ \psi_{\mt W.,j;(\mt D.')} \right] \; =\;
\left[ \psi_{\mt V.,i;(\mt D.)}^* \ , \ \psi_{\mt W.,j;(\mt D.')} \right]
               &=& 0\;,\qquad
V \neq W\;,\; V,W \in \ot G. \ .  \label{31.4}
\end{eqnarray}

\begin{rem}\label{not-tensor}
  For every $V \in \ot G.$ put $d(V):=\mathrm{dim}\al H._\mt V.\in\N$
  and denote by $\al O._{d(V)}$ the Cuntz algebra generated by $\al
  H._\mt V.$, and by $^0 \al O._\mt {d(V)}.$ the dense *-subalgebra of
  $\al O._\mt {d(V)}.$ algebraically generated by $\al H._\mt V.$
  (cf.~Section~\ref{Cuntz_Pimsner}).  By universality of the Cuntz
  relations, we find that for every $V \in \ot G.$ there is a unital
  *-monomorphism $^0 \al O._\mt {d(V)}. \hookrightarrow \ ^0 \al F.$.
  Moreover, (\ref{31.4}) implies that the algebraic tensor product $\
  ^0 \al O._{\ot G.} :=$ $\odot_\mt V. \ ^0 \al O._\mt {d(V)}.$ is
  embedded in $\ ^0 \al F.$. We denote by
\begin{equation}
\label{def_otg}
\al O._{\ot G.} := \mathop{\bigotimes}\limits_{\mt{V \in \ot G.}.} \al O._{d(V)}
\end{equation}
the {\it C*}-tensor product of the corresponding
Cuntz algebras (which are nuclear).
From Eq.~(\ref{31.3}) we have that every Hilbert $\ot R.$-bimodule
$\ot H._\mt V.$ defined in Section~\ref{SecLie} is also embedded in $\ ^0 \al F.$.
Recall from Remark~\ref{GenCP} that
$^0 \al O._{\ot H._V}$ is the dense *-subalgebra
of $\al O._{\ot H._V}$ algebraically generated by $\ot H._\mt V.$.
Then, again by universality, we have that $^0 \al O._{\ot H._V}$ is embedded
in $^0 \al F.$.
Note, nevertheless, that in general the *-algebras
$^0 \al O._{\ot H._{V_1}}$, $^0 \al O._{\ot H._{V_2}}$, $V_1,V_2\in\ot G.$,
do {\em not} appear in $\ ^0 \al F.$ as a tensor product
due to the twist introduced by the chain group action
in Eq.~(\ref{31.3}). In fact, it is
easy to check using the previous relations,
that if the chain group action described by $\alpha$ is
nontrivial, then
\[
\psi_2 (B \psi_1) = \alpha(B) \psi_1 \psi_2\not= (B \psi_1) \psi_2\,,
\quad\psi_i\in\al O._{\ot H._{V_i}},\, i=1,2;\, B\in\ot R.\,.
\]
\end{rem}

In the following we will show that $^0 \al F.$ admits a nontrivial {\it C*}-norm.
Let
\[
\T^\infty:=\times_{\mt V \in {\ot G.}.} \T
\ \ , \ \
\Z^\infty := \times_{\mt V \in {\ot G.}.} \Z
\]
be the product of circles $\T$ (resp. $\Z$) and
denote its elements as maps $z\colon\ot G. \rightarrow \T$
(resp. $k\colon\ot G. \rightarrow \Z$).
On $^0 \al F.$ we have the following natural actions by *-automorphisms
\begin{eqnarray}
\label{eq_g}
\beta \colon \al G. \rightarrow \aut ^0 \al F.
\ \ , \ \
\beta_g ( \psi_{\mt V.,i;(\mt D.)} ) := V(g) \psi_{\mt V.,i;(\mt D.)}
                                    = V_\mt D.(g) \psi_{\mt V.,i;(\mt D.)}
     \quad\mathrm{and}\quad
      \beta_g (B):=B\\
\label{eq_t}
\delta\colon \T^\infty \rightarrow \aut ^0 \al F.
\ \ , \ \
\delta_z ( \psi_{\mt V.,i;(\mt D.)} ) :=  z(V) \cdot \psi_{\mt V.,i;(\mt D.)}
     \quad\mathrm{and}\quad
     \delta_z (B):=B\,,\;B\in\ot R.\, .
\end{eqnarray}
Moreover, we consider the sets
\begin{eqnarray*}
\Z_0^\infty &:=&\{k \in \Z^\infty \mid \mr supp.(k)<\infty\} \\
\N_0^\infty &:=&\{r\in\Z_0^\infty\mid  r(V)\geq 0\,,V\in\ot G.\}\,.
\end{eqnarray*}
For every $k \in \Z_0^\infty$, we define the spectral subspaces
w.r.t.~the action $\delta$ as follows:
\[
^0 \al F.^k
:=
\left\{
T \in \ ^0 \al F. \mid
\delta_z (T) = \prod_{V \in \ {supp} (k)} z(V)^{k(V)} \ T
\;,\;\; z\in\T^\infty
\right\} \ .
\]
Note that these subspaces are a natural generalization of the standard spectral
subspace considered for a single Cuntz-Pimsner algebra in (\ref{SpecS1}).
If $k_1\not=k_2$, then the
corresponding spectral subspaces have a trivial intersection.
Moreover, if $0 \in \Z_0^\infty$ denotes the zero section,
then $^0\al F.^0$ is the fixed-point
algebra w.r.t. the $\T^\infty$-action.
It is natural to introduce the notation
$\ ^0 \al O.^k :=\,{^0 \al F.}^k \cap\,  {^0\al O.}$, $k \in \Z_0^\infty$,
where ${^0\al O.}$ is the algebraic tensor product of the algebraic
part of the Cuntz algebras
defined in Remark~\ref{not-tensor}.
We have $\ ^0 \al F. =\mr span.\{ ^0 \al F.^k\mid k \in \Z_0^\infty\}$
as well as $\ ^0 \al O._{\ot G.} =\mr span. \{^0 \al O.^k\mid k \in \Z_0^\infty\}$.
Note that for every $V \in \ot G.$,
it turns out that $\ ^0 \al O._{d(V)} \cap$ $^0 \al F.^k$ coincides with the
spectral subspace $\ ^0 \al O._{\mt d(V).}^{\mt k(V).}$ w.r.t. the standard
canonical circle action on $\ ^0 \al O._{d(V)}$.
We introduce the projection
\begin{equation}
\label{eq_m0}
m_0 : \ ^0 \al F. \rightarrow \ ^0 \al F.^0
\ \ , \ \
m_0 \left( \sum^{\mbox{\tiny fin}}_{\mt k \in \Z_0^\infty.} T_k \right) := T_0 \,,
\end{equation}
where $T_k$ is in the corresponding spectral subspace $^0 \al F.^k$, and
the symbol $\sum^{\mbox{\tiny fin}}$ indicates a sum over a finite set
of indexes.
If $\al V.$ is a finite subset of $\ot G.$ and
$T \in$ $^0 \al F.^k$ with $k(V) \neq 0$ for some
$V \in \al V.$, then $m_0 (T) = 0$.

For any $r \in \N_0^\infty$ put $\al V. :=$ ${supp}(r)$ which is a finite subset
of $\ot G.$. Then
$\al H.^r_{\mt{\mathcal{V}}.}$ denotes the linear span of elements of
the form $\prod_{\mt{V \in \al V.}.} \varphi_\mt V.$, with
$\varphi_\mt V. \in\al H._\mt V.^{\mt r(V).}$.
It is clear that we may identify $\al H.^r_{\mt{\mathcal{V}}.}$ with the
tensor product $\otimes_{\mt{V \in \al V.}.} \al H._{\mt V.}^{\mt r(V).}$
(recall that $\al H._{\mt V.}^{\mt r(V) .}$ is identified with the tensor product
$\otimes^{\mt r(V) .} \al H._{\mt V.}$).
We denote by
$\al L.(\al H._{\mt{\al V.}.}^r , \al H._{\mt{\al V.}.}^s )$
the space of linear maps from
$\al H._{\mt{\al V.}.}^r$ into $\al H._{\mt{\al V.}.}^s$:
it coincides with the linear span of
elements of the form $\varphi' \varphi^*$, where $\varphi' \in \al H._{\mt{\al V.}.}^s$,
$\varphi \in \al H._{\mt{\al V.}.}^r$, where
$r,s \in \N_0^\infty$ have support $\al V.$.

\begin{rem}
\label{rem_part1_0}
For every finite $\al V. \subset \ot G.$, $k \in \Z_0^\infty$,
we consider the vector space
\[
\ot H._{\al V.}^k
:=
{\mathrm{span}}
\prod_{V \in \al V.} \ot H._\mt V.^{k(V)}
\]
(for $k(V) < 0$, $V \in \al V.$, we define as usual
$\ot H._\mt V.^{k(V)} :=$ $\left( \ot H._\mt V.^{-k(V)} \right)^*$
which can be interpreted as the dual bimodule of bounded, right $\ot R.$-module
maps from $\ot H._V^{-k(V)}$ into $\ot R.$;
for $k(V) = 0$, recall that by definition $\ot H._V^0 = \ot R.$).
Now, we may identify $\ot H._{\al V.}^k$ with the $\ot R.$-bimodule
inner tensor product
\[
\bigotimes_{\mt{V \in \al V.}.} \ot H._\mt V.^{k(V)} \ \
\]
\end{rem}

In the following lemma we give several useful characterizations of
the spectral subspaces $^0 \al F.^k$. These will play a fundamental
role in the rest of this section.

\begin{lem}
\label{lem_fk}
For every $k \in \Z_0^\infty$, we have
\begin{eqnarray}
\label{eq_dec}
^0 \al F.^k &=&
\left\{  \sum_{\mt{\al V.}.}^{\mbox{\tiny fin}}  B_{\mt{\al V.}.} T_{\mt{\al V.}.}
          \mid
           B_{\mt{\al V.}.} \in \ot R.
           \ , \
           T_{\mt{\al V.}.} \in
             \al L.( \al H._{\mt{\al V.}.}^r , \al H._{\mt{\al V.}.}^{r+k} )
           \;,\;\; r\in\N_0^\infty \ , \ \al V. \subseteq \gfin
           \ {\mr finite.} \
\right\}\\
\label{eq_dec1}
            &=&
\left\{  \sum_{\mt{\al V.}.}^{\mbox{\tiny fin}}  R_{\mt{\al V.}.}
         \mid
          R_{\mt{\al V.}.} \in
             \al L.( \ot H._{\mt{\al V.}.}^r , \ot H._{\mt{\al V.}.}^{r+k} )
           \;,\;\; r\in\N_0^\infty \ , \ \al V. \subseteq \gfin
           \ {\mr finite.} \
\right\}\ .
\end{eqnarray}
In particular, with the conventions of Remark~\ref{rem_part1_0}
we have $\ot H._{\al V.}^k \subset\, {^0\al F.^k}$. Moreover
recalling Remark~\ref{GenCP} we can also write
\[
^0 \al F.^k = \mr span.\left(\prod_{V \in \ {supp} (k)}
              { ^0 \al O._{\ot H._\mt V.}^{k(V)} }
              \right)\,.
\]
\end{lem}
\begin{beweis}
From the structure of the spectral subspaces of the single
Cuntz-Pimsner algebras (cf.~Remark~\ref{GenCP}) it is clear
that every element of $^0 \al F.^k$ is the sum of terms of the form
\[
T = B_0 T_{\mt V_1 .} B_1 \cdots B_{n-1} T_{\mt V_n .} B_n\,,
\]
where $V_1,\ldots, V_n \in \ot G.$,
$T_{\mt V_i .} \in \al L.( \al H._{\mt V_i.}^{\mt r(V_i) .} ,
\al H._{\mt V_i.}^{\mt s(V_i) .} )$,
$s(V_i)=r(V_i) +k(V_i)$, and $B_i \in \ot R.$.
By using Eq.~(\ref{31.4}), we may change the order of terms
$T_{\mt V_i .}$, $T_{\mt V_j .}$, if $i\neq j$, and put together the
terms arising from the same representation.
Thus, without loss of generality, we may assume that $V_i \neq V_j$
if $i \neq j$.
Moreover, every $T_{\mt V_i .}$ can be written as a sum of terms of the form
$\psi_{\mt V_i.} \varphi^*_{\mt V_i.}$, where
$\psi_{\mt V_i.} \in\al H._{\mt V_i.}^{\mt s(V_i) .}$,
$\varphi_{\mt V_i.} \in \al H._{\mt V_i.}^{\mt r(V_i) .}$.
Therefore applying Eq.~(\ref{31.3}) successively we can write all elements
from $\ot R.$ to the left, i.e.~for some $B \in \ot R.$ we have
\[
T = BT'
\quad\mr with.\quad
T' := T_{\mt V_1 .} \cdots T_{\mt V_n .} \ \, .
\]
Finally, (\ref{31.4}) implies that $T'$
can be rewritten as
\[
\psi_{\mt V_1.} \psi_{\mt V_2.} \cdots \psi_{\mt V_n.}
\cdot
\varphi^*_{\mt V_1.} \varphi^*_{\mt V_2.} \cdots \varphi^*_{\mt V_n.}
\ .
\]
Hence $T' \in\al L.( \al H._{\mt{\al V.}.}^r , \al H._{\mt{\al V.}.}^{r+k} )$,
where
$\al V. :=$
$\left\{  V_1 , \ldots , V_n  \right\}$.

The second characterization of the spectral subspace $^0 \al F.^k$
follows from the fact that
\[
 B T_{\mt{\al V.}.} \in\al L.( \ot H._{\mt{\al V.}.}^r , \ot H._{\mt{\al V.}.}^{r+k} )
\]
for every
$B \in \ot R.$, $T \in\al L.( \al H._{\mt{\al V.}.}^r , \al H._{\mt{\al V.}.}^{r+k} )$.

The last statements in the lemma follow immediately from the structure of the
spectral subspaces of the single Cuntz-Pimsner algebras
(cf.~Remark~\ref{GenCP}).
\end{beweis}

\begin{defi}
\label{def_suppt}
For every $T \in$ $^0 \al F.^k$, $k \in \Z_0^\infty$, we denote by
${supp} (T)$ the finite subset in $\ot G.$ given by the union of all
$\al V.\subset\ot G.$ that appear in the decomposition (\ref{eq_dec})
or (\ref{eq_dec1}).
\end{defi}

\subsubsection{Crossed products}

In the following we will study Hilbert space representations of the
*-algebra $^0\al F.$. We begin introducing the notion of covariant
representation associated to a family of unitary finite-dimensional
representations of the compact group $\al G.$ and a fixed chain group
action on a {\it C*}-algebra.

\begin{defi}
\label{def_cv}
Let $\ot G.=\ot G.(\al G.)$ be a family of finite-dimensional unitary
representations of $\al G.$ and denote by $\al O._{\ot G.}$ the tensor
product of the corresponding Cuntz algebras (cf.~(\ref{def_otg})).
Consider a unital {\it C*}-algebra $\ot R.$ and a fixed
chain group action $\alpha\colon\ot C.(\al G.) \rightarrow\aut \ot R.$.
A {\bf covariant representation} of $(\ot R.,\ot G.,\ot C.(\al G.),\alpha)$ is
a triple $(\ot h.,\pi,\eta)$ where
$\pi\colon \ot R. \to \al L.(\ot h.)$ and
$\eta\colon \al O._{\ot G.} \to \al L.(\ot h.)$ are non-degenerate
Hilbert space representations satisfying
\[
\eta ( \psi_{\mt V.,i;(\mt D.)} ) \cdot \pi(B)
   = \pi \circ \alpha_{\mt [D].} (B)
     \cdot \eta ( \psi_{\mt V.,i;(\mt D.)} )\;,\quad
     B \in \ot R.\;,[D] \in \ot C.(\al G.)\;,
     \psi_{\mt V.,i;(\mt D.)} \in\al H._{\mt V.} \subset\al O._{\ot G.}\,.
\]
\end{defi}

Here we identify $\psi_{\mt V.,i;(\mt D.)}$
with the corresponding image in $\al O._{\ot G.}$ according to the natural
inclusions $\al H._\mt V.
\hookrightarrow\al O._{\mt {d(V)}.} \hookrightarrow \al O._{\ot G.}$. Now, since
$\al O._{\ot G.}$ is a tensor product
labeled by the different elements in $\ot G.$
and since $\eta$ is non-degenerate
(i.e.~unital), we find that the operators
$\eta ( \psi_{\mt V.,i;(\mt D.)} )$ satisfy the relations
(\ref{31.1}), (\ref{31.2}) and (\ref{31.4}).
Moreover from Definition~\ref{def_cv} it is also clear that
the operators $\pi(B)$, $\eta (\psi_{\mt V.,i;(\mt D.)})$
satisfy the relation (\ref{31.3}).
Therefore, defining
\begin{equation}
\label{eq_ir}
\pi \rtimes \eta\colon \ ^0 \al F. \to \al L.(\ot h.)
\ \ , \ \
\pi \rtimes \eta (B \psi_{\mt V.,i;(\mt D.)} ) :=
\pi(B) \cdot \eta ( \psi_{\mt V.,i;(\mt D.)} )
\ ,
\end{equation}
we obtain a representation of $^0 \al F.$ on the Hilbert space
$\ot h.$. Let us define the
{\it C*}-algebra $\al F.$ as the closure of $^0 \al F.$ w.r.t. the
{\it C*}-norm
\begin{equation}\label{CNorm}
  \|F\|:=
  \mathop{\mr sup.}\limits_{(\pi,\eta)}\|(\pi \rtimes \eta)(F) \|_{\mrt op.}
          \;, F\in \ ^0 \al F.\;.
\end{equation}

Thus by construction we have that for every covariant
representation $(\pi,\eta)$, there exists a representation
$\Pi$ of $\al F.$ extending the former pair, i.e.
\begin{equation}
\label{eq_ir1}
\Pi\colon \al F. \to \al L.(\ot h.)
\quad\mr with.\quad
\Pi (T) = \eta (T)
\ , \
\Pi(B) = \pi(B) \;,\quad
T \in \ ^0 \al O._{\ot G.} \subset\ ^0 \al F.\;,\;B \in \ot R.\,.
\end{equation}

The proof of the following theorem uses the notion of twisted tensor product
introduced by Cuntz in \cite[\S 1]{Cun81}.

\begin{teo}
\label{thm_exf}
Let $\al G.$ be a compact group, $\ot R.$ a unital {\it C*}-algebra
and $\alpha\colon\ot C.(\al G.) \rightarrow\aut \ot R.$ a fixed
chain group action.
Then, for every set $\ot G.$ of finite-dimensional representations of $\al G.$,
there exists a universal {\it C*}-algebra
\[
\ot R. \rtimes^\alpha \ot G.
\]
generated by $\ot R.$ and $\{\psi_{\mt V.,i;(\mt D.)}\}_{\mt {V\in\ot G.}.}$
and satisfying the relations (\ref{31.1})-(\ref{31.4})
and such that (\ref{eq_g}),(\ref{eq_t}) extend to automorphic actions on
$\al F.$.
\end{teo}
\begin{beweis}
We will divide the proof into several steps:

1. Consider first a covariant representation
   $(\ot h., \pi_0 , U )$ of the {\it C*}-dynamical system
   $(\ot R.,\ot C.(\al G.),\alpha)$
   in the usual sense of crossed products by group actions,
   i.e. $\pi_0\colon \ot R. \to \al L.(\ot h.)$ is a non-degenerate representation
   and $U\colon \ot C. (\al G.) \to \al U.(\ot h.)$ is a unitary representation
   of the chain group satisfying
\[
U_{\mt [D].} \cdot \pi_0 (B) = \pi_0 \circ \alpha_{\mt [D].} (B) \cdot U_{\mt [D].}
\;,\quad B \in \ot R.\;,\; [D] \in \ot C. ( \al G. )\;.
\]
   It is well-known that
   $\left\| B \right\| =\mathop{\sup}\limits_{\pi_0} \left\| \pi_0(B) \right\|$
   (see \cite{bPedersen79}). Moreover, consider a Hilbert space representation
   $\eta_0\colon \al O._{\ot G.} \to \al L.(\ot h._0)$. By simplicity of the Cuntz
   algebras, we have that $\al O._{\ot G.}$ is also simple, hence
   $\eta_0$ is faithful.

2. In terms of $\eta_0$ and the previous covariant representation
   $\pi_0$ we will introduce next a covariant representation of $(\ot
   R.,\ot G.,\ot C.(\al G.),\alpha)$ in the sense of
   Definition~\ref{def_cv}. We define
\[
\pi \colon \ot R. \to \al L.( \ot h. \otimes \ot h._0 )
\ \ , \ \
\pi (B) := \pi_0 (B) \otimes \1
\ ,
\]
\[
\eta \colon \al O._{\ot G.} \to \al L.( \ot h. \otimes \ot h._0 )
\ \ , \ \
\eta ( \psi_{\mt V.,i;(\mt D.)} ) :=
U_{\mt [D].} \otimes \eta_0 ( \psi_{\mt V.,i;(\mt D.)} )
\ .
\]
Note that with this definition we have
\begin{eqnarray*}
\eta(\psi_{\mt V.,i;(\mt D.)})\cdot\pi(B)
    &=& \Big( U_{\mt [D].} \otimes \eta_0 ( \psi_{\mt V.,i;(\mt D.)} ) \Big)
         \cdot \Big( \pi_0 (B) \otimes \1 \Big) \\
    &=& \Big( U_{\mt [D].} \cdot \pi_0 (B) \Big)
        \otimes\eta_0 ( \psi_{\mt V.,i;(\mt D.)} ) \\
    &=& \Big( \pi_0 \circ \alpha_{\mt [D].} (B) \cdot U_{\mt [D].} \Big)
        \otimes\eta_0 ( \psi_{\mt V.,i;(\mt D.)} )\\
    &=& \Big(\pi_0 \circ \alpha_{\mt [D].} (B) \otimes \1 \Big)
        \cdot \Big( U_{\mt [D].} \otimes \eta_0 ( \psi_{\mt V.,i;(\mt D.)} ) \Big)\\
    &=& \Big(\pi\circ\alpha_{\mt [D].}\Big)(B)\cdot \eta(\psi_{\mt V.,i;(\mt D.)}) \,,
\end{eqnarray*}
This shows that
$(\pi,\eta)$ specifies a representation of the relation (\ref{31.3}).
Furthermore we verify that $\eta$ is also a representation of the remaining
relations:
\[
\sum_i
\eta ( \psi_{\mt V.,i;(\mt D.)} ) \eta ( \psi_{\mt V.,i;(\mt D.)} )^* =
\sum_i
\Big( U_{\mt [D].} U_{\mt [D].}^* \Big)
\otimes
\eta_0 ( \psi_{\mt V.,i;(\mt D.)}  \psi_{\mt V.,i;(\mt D.)}^* ) =
\1 \otimes \1 \ ,
\]
\[
\eta ( \psi_{\mt V.,i;(\mt D.)} )^* \eta ( \psi_{\mt V.,j;(\mt D'.)} ) =
\Big( U_{\mt [D].}^*  U_{\mt [D'].} \Big) \otimes
\eta_0 ( \psi_{\mt V.,i;(\mt D.)}^* \psi_{\mt V.,j;(\mt D'.)}  ) =
\delta_{ij}\; (\1 \otimes \1) \ ,
\]
and
\begin{eqnarray*}
\eta ( \psi_{\mt V.,i;(\mt D.)} ) \eta ( \psi_{\mt W.,j;(\mt D'.)} )
  & = & \Big(U_{\mt [D].} \otimes \eta_0 ( \psi_{\mt V.,i;(\mt D.)} ) \Big)
        \cdot \Big(U_{\mt [D'].} \otimes \eta_0 ( \psi_{\mt W.,j;(\mt D'.)} ) \Big)  \\
  & = & \Big( U_{\mt [D].}  U_{\mt [D'].} \Big)
        \otimes \eta_0 (  \psi_{\mt V.,i;(\mt D.)}   \psi_{\mt W.,j;(\mt D'.)}   ) \\
  & = & \Big( U_{\mt [D'].} U_{\mt [D].} \Big) \otimes
        \eta_0 (  \psi_{\mt W.,j;(\mt D'.)}  \psi_{\mt V.,i;(\mt D.)}   ) \\
  & = & \eta ( \psi_{\mt W.,j;(\mt D'.)} ) \eta ( \psi_{\mt V.,i;(\mt D.)} ) \ ,
\end{eqnarray*}
where we used (\ref{31.4}) and the fact that $\ot C. (\al G.)$ is
Abelian. In the same way, one can verify the relations
\[
 \eta ( \psi_{\mt W.,j;(\mt D'.)} )^* \eta ( \psi_{\mt V.,i;(\mt D.)} )
   =\eta ( \psi_{\mt V.,i;(\mt D.)} ) \eta ( \psi_{\mt W.,j;(\mt D'.)})^*\,.
\]

3. In the preceding step we have shown that
for any covariant representation $(\ot h.,\pi_0,U)$
and any representation $\eta_0\colon\al O._\ot G.\to\al L.(\ot h._0)$
as in step~1 we have constructed a covariant representation $(\pi,\eta)$ of
$(\ot R.,\ot G.,\ot C.(\al G.),\alpha)$ over the Hilbert space
$\ot h.\otimes\ot h._0$. In particular, $(\pi,\eta)$ defines also a
representation $\pi \rtimes \eta$ of
$^0\al F.$ on $\ot h.\otimes\ot h._0$. Therefore,
we can introduce the {\it C*}-norm as in Eq.~(\ref{CNorm})
\[
\|F\|:=\mathop{\mr sup.}\limits_{(\pi,\eta)}\|(\pi \rtimes \eta)(F) \|_{\mrt op.}
          \;, F\in \ ^0 \al F.\;,
\]
and complete $^0 \al F.$ w.r.t.~this norm:
\[
\ot R. \rtimes^\alpha \ot G.:=\clo_{\|\cdot\|} \ \Big(\ ^0\al F.\Big)\,.
\]
(We will sometimes denote $\ot R. \rtimes^\alpha \ot G.$ simply by
$\al F.$.)

Note that the preceding {\it C*}-norm extends the {\it C*}-norm
on $\ot R.$, since recall that we have
\[
\left\| B \right\|= \sup_{\pi_0} \left\| \pi_0 (B) \right\|_{\mrt op.}
                  =\sup_\pi \left\| \pi (B)\right\|_{\mrt op.}
                  =\sup_{(\pi,\eta)}\|(\pi \rtimes \eta)(B) \|_{\mrt op.}\,.
\]

To show the universal property, let $\al F.'$ be a {\it C*}-algebra with
generators satisfying the relations (\ref{31.1})-(\ref{31.4}), and $\Pi'
\colon \al F.' \to \al L.(\ot h.)$ a faithful, non-degenerate representation.
Then, $\ot R.$ and $\al O._{\ot G.}$ are {\it C*}-subalgebras of $\al
F.'$ and
$( \pi' := \Pi' \rest_{\ot R.} \ , \ \eta' := \Pi' \rest_{\al O._{\ot G.}} )$
specifies a covariant representation. Since
\[
 \|(\pi' \rtimes \eta')(F) \|_{\mrt op.}\leq   \|F\|\;,\quad F\in \ ^0\al F.
\]
it is clear that there is a monomorphism $\al F.\hookrightarrow \al F.'$\,.

4. Finally we address the question of
the automorphic extensions. Let $g \in \al G.$ and
$\beta_g$ defined as in (\ref{eq_g}) in terms of
unitary representations. Note that if
$( \pi , \eta )$ is a covariant representation, then
$( \pi , \eta \circ \beta_g )$ is again a covariant
representation. From this and the universality of the
Cuntz algebra it is clear that $\beta$ extends to an
automorphic action of the compact group $\al G.$.
A similar argument shows that (\ref{eq_t}) extends to an
automorphic action of $\T^\infty$.
\end{beweis}

\begin{rem}
\label{ex_rg}
\begin{itemize}
\item[(i)] The {\it C*}-algebra $\ot R. \rtimes^\alpha \ot G.$
generalizes several well-known constructions.
If $\ot G.$ has a unique element $V$, then $\ot R.
\rtimes^\alpha \ot G.$ coincides with the Cuntz-Pimsner algebra
$\al O._{\ot H._V}$ studied in Section~\ref{SecLie}.
In particular, if $V = D$
is an irreducible representation, then $\ot R. \rtimes^\alpha \ot G.$
is isomorphic to the Stacey crossed product $\ot R. \rtimes^{d(\mt
  D.)}_\alpha \N$ (see \cite{Stacey}).  If $\al G.$ is Abelian and
$\ot G. = \ot C. (\al G.)$, then $\ot R. \rtimes^\alpha \ot G.$ is
isomorphic to the usual crossed-product $\ot R. \rtimes_\alpha \ot
C.(\al G.)$ by the $\ot C.(\al G.)$-action on $\ot R.$.
Finally, if the action of the chain group on $\ot R.$ is trivial,
then $\ot R. \rtimes^\alpha \ot G.$ reduces to the
tensor product $\ot R. \otimes \al O._{\ot G.}$ (recall that $\al
O._{\ot G.}$ is nuclear, thus we do not have to specify the norm of
the tensor product).

\item[(ii)] Note also that the {\it C*}-algebra $\ot R. \rtimes^\alpha \ot G.$
 is generated as a {\it C*}-algebra by the family of Cuntz-Pimsner algebras
 \[
   \{\al O._{\ot H._\mt V.}\}_{\mt {V\in\ot G.}.}\,,
 \]
 where the product of the elements in different algebras is twisted by
 the chain group action. The present construction has also some elements
 of Cuntz notion of twisted tensor product given in \cite{Cun81}.
 Finally, note also that
 the {\it C*}-norm restricted to a single Cuntz-Pimsner algebra is unique.
\end{itemize}
\end{rem}

To simplify notation we write $\al F. := \ot R.
\rtimes^\alpha \ot G.$.  We denote by $\al F.^k$, $k \in \Z_0^\infty$,
the spectral subspaces of $\al F.$ w.r.t.~the $\delta$-action. It is
clear that each $\ ^0 \al F.^k$ is a dense subspace of $\al F.^k$.
Moreover, the projection onto the zero spectral subspace in
(\ref{eq_m0}) extends in a natural way to $m_0\colon\al F. \to\al F.^0$.

By construction, for every $V \in \ot G.$ there corresponds an
algebraic Hilbert space $\al H._\mt V.$ with support $\1$ contained in
$\al F.$ (see Eq.~(\ref{31.1})).  Moreover, from (\ref{eq_g}) we have
that $\beta_g(\psi) =V(g) \psi$, $g \in \al G.$, $\psi \in \al H._\mt
V.$, and for every $V \in \ot G.$, there is an endomorphism
$\sigma_\mt V. \in \mr End. \al F.$ induced as usual by
\[
\sigma_\mt V. (T) :=
\sum_i \psi_{\mt V.,i;(\mt D.)} T \psi^*_{\mt V.,i;(\mt D.)}
\ \ , \ \
T \in \al F.
\ ,
\]
where $\left\{ \psi_{\mt V.,i;(\mt D.)} \right\}_i$ is an adapted orthonormal
basis of the algebraic Hilbert space $\al H._\mt V.$.
Each $\sigma_\mt V.$ coincides with the canonical endomorphism
on the Cuntz algebra $\al O._{d(V)} \subset\al F.$.
Let $\al A. \subset \al F.$ denote the fixed-point algebra w.r.t. the action
$\beta$. Since the construction of $\sigma_\mt V.$ does not depend on the choice of
the basis $\psi_\mt V.$, we obtain
$\sigma_\mt V. \circ \beta_g =\beta_g \circ \sigma_\mt V.$,
$g \in \al G.$, as in Section~\ref{SecLie}.
This implies that every $\sigma_\mt V.$ restricts to a
canonical endomorphism $\rho_\mt V. \in \mr End.\al A.$.

\subsection{Minimal and regular {\it C*}-dynamical systems}
\label{5.2}

We now assume that $\ot R.$ is Abelian and denote
it by $\al Z.$ as in the first part of the paper.
Let us consider the class $\gfin$ of finite-dimensional representations
introduced in Definition~\ref{def_gfin}.
For the choice $\ot G. :=\gfin$ we obtain from the construction in the
preceding subsection a {\it C*}-dynamical system $( \al F. ,\al G.,\beta)$,
$\al F. :=$ $\al Z. \rtimes^\alpha \gfin$,
with $\T^\infty$-action given by $\delta$ (cf.~Eq.~(\ref{eq_t})) and
fixed-point algebra $\al A.$. Recall also
that $\al F.$ contains
all Cuntz-Pimsner algebras $\al O._{\ot H._\mt V.}$, $V \in\gfin$.
Moreover, by definition of the group action given in Eq.~(\ref{eq_g}) it is clear
that $\beta$ restricts to each single Cuntz-Pimsner algebra, so that the large
{\it C*}-dynamical system $( \al F. ,\al G.,\beta)$ contains all single
{\it C*}-dynamical systems $(\al O._{\ot H._\mt V.},\al G.)$
with corresponding fixed-point algebra $\al A._\mt V.$ considered
in Theorem~\ref{thm_main}. Note also that we have the following relation
between the fixed-point algebras of the single systems and the
fixed-point algebra $\al A.$ of the large
{\it C*}-dynamical system:
\[
 \al A._\mt V.\subset \al A. \;,\quad V\in\gfin\,.
\]

In the present subsection the relative commutant $\al A.'\cap\al F.$
will analyzed in detail. We begin stating the following stability
result:
\begin{lem}\label{StableRel}
The relative commutant $\al C.:=\al A.'\cap\al F.$ of the
{\it C*}-dynamical system $(\al F. ,\al G.,\beta)$ with
fixed-point algebra $\al A.$ is stable under the
$\T^\infty$-action given by $\delta$. Therefore
$\al C.$ is generated as a {\it C*}-algebra by the spectral
subspaces
\[
 \al C.^k:=\al F.^k\cap\al C.\;,\quad k\in\Z^\infty_0.
\]
\end{lem}
\begin{beweis}
Let $z \in {\mathbb{T}}^\infty$ and $C \in \al C.$.
Since $\al A.$ is generated by the corresponding spectral
subspaces, it is enough to show that
\[
 \delta_z(C) A = A \delta_z(C)\quad\mr for~any.\quad
                 A\in\al A.^k\;,\;\; k\in\Z^\infty_0\,.
\]
But this follows from
\[
\delta_z (C) \prod_\mt V. z(V)^\mt {n(V)}. A =
\delta_z (C A) = \delta_z (A C)
 = \prod_\mt V. z(V)^\mt {n(V)}. A \delta_z(C)\,.
\]
This shows the stability of $\al C.$ and the fact that it is generated by the spectral
subspaces $\al C.^k$.
\end{beweis}

In the next proposition, we show that {\em any} finite-dimensional
representation of $\al G.$  is realized on an algebraic
Hilbert space in {\it C*}-algebra $\al F.$.

\begin{pro}
Let $V$ be a finite-dimensional representation of $\al G.$
(not necessarily belonging to $\gfin$) on $\al H._\mt V.$.
Then, $\al H._\mt V.$ can be identified with an algebraic Hilbert space
in $\al F.$ with support $E_\mt V.$. Moreover,
$\beta_g (\psi) = V(g) \psi$, $g \in \al G.$,
$\psi \in \al H._\mt V.$.
\end{pro}
\begin{beweis}
Let $V_\mt D.$ be an irreducible subrepresentation of $V$ of class $D$.
Define $L := V \oplus V_\mt D.$
and consider the representation
$W := L \oplus \overline{\mathrm{det}L}$. By construction
$W \in \gfin$ and the algebraic Hilbert space
$\al H._\mt W.$ is contained in $\al F.$ with support $\1$.
(Note that $W \in \gfin$ also in the case where $V$ is an irreducible
representation of an Abelian group.)
There is a reducing projection $E_\mt V.$ associated with the decomposition
$\al H._\mt W. = \al H._\mt V.
\oplus \al H._\mt {V_\mt D. \oplus \overline{\mathrm{det}L}}.$.
Thus, $\al H._\mt V.$ has support $E_\mt V.$ in $\al F.$.
Moreover,
$\beta_g (\psi) =
 ( V(g) \oplus V_\mt D.(g) \oplus \overline{\mathrm{det}L(g)} ) (E_\mt V. \psi) =
 V(g) \psi$, $g \in \al G.$,
$\psi \in \al H._\mt V.$.
\end{beweis}

From Definition~\ref{def_gfin} for every $V \in \gfin$
there exists $n(V) \in {\mathbb{N}}$ and a $\al G.$-invariant
isometry $S_\mt V.$ of $\al H._\mt V.^\mt {n(V)}.$.
In the following lemma we will summarize some useful properties
of these invariant vectors (see also Lemma~\ref{Decomp1}).

\begin{lem}\label{InvVect}
For every $V \in \gfin$, there
exists $S_\mt V. \in \al H._\mt V.^\mt {n(V)}. \subset \al F.$ such that
\begin{itemize}
\item[(i)]
$S_\mt V.^* S_\mt V. = \1$, $\beta_g (S_\mt V.) = S_\mt V.$
and $S_\mt V. Z = Z S_\mt V.$,
$Z \in \al Z. \ , \ g \in \al G.$.
\end{itemize}
Moreover,
\begin{itemize}
\item[(ii)]
if $V,W\in \gfin$, $V \neq W$, then $S_\mt V. T = T S_\mt V.$,
$T \in \al O._{\ot H._\mt W.} \subset \al F.$.
\item[(iii)]
if $\delta$ is the $\T^\infty$-action on $\al F.$ (cf.~Eq.(\ref{eq_t}))
and $C\in\al C.:=\al A.'\cap\al F.$, then
$\delta_z (C) S_\mt V. = S_\mt V. \delta_z (C)$,
$z\in\T^\infty$.
\end{itemize}
\end{lem}
\begin{beweis}
The existence of the elements $S_\mt V.$ and first two equations in (i)
follow from Definition~\ref{def_gfin}~(ii)
and Lemma~\ref{Decomp1}. Therefore $S_\mt V.\in\al A._\mt V.\subset\al A.$, where
$\al A._\mt V.$ is the fixed-point algebra of the {\it C*}-dynamical system
$(\al O._{\ot H._\mt V.},\al G.)$. The last equation in
(i) follows from the equation $\al Z.=\al A._\mt V.'\cap\al A._\mt V.$
(cf.~Theorem~\ref{thm_main}~(ii)).
Recall that the Cuntz-Pimsner algebra $\al O._{\ot H._\mt W.}$ is generated by
$\al Z.$ and $\al H._\mt W.$. Therefore the equation in (ii) follows from
$S_\mt V.\in\al Z.'\cap\al A.$ (cf.~part~(i)) and relation (\ref{31.4}).
Part (iii) is an immediate consequence of Lemma~\ref{StableRel} and
the fact that $S_\mt V.\in\al A.$.
\end{beweis}

For every finite set $\al W. \subseteq \gfin$, we define the isometry
\begin{equation}
\label{eq_SW}
S_{\al W.}
:=
\prod_{\mt{ W \in \al W.}.} S_W
\end{equation}
and for every $T \in \al F.$, $p \in \N$, we put
\begin{equation}
\label{def_ttilde}
\widetilde T_{\al W.,p} :=
(S_{\al W.}^p)^* \ T \ S_{ \mt{\al W.}.}^p
\ \in \al F. \ .
\end{equation}

\begin{lem}
\label{lem_rv}
Let $R \in$ $^0 \al F.^k$, $k \in \Z_0^\infty$, and set
$\al V. :={supp} (R)$ (see Definition~\ref{def_suppt}).
Then, for every finite
set $\al W. \subseteq \gfin$ such that $\al V. \subseteq \al W.$,
and for every $p \in \N$, we have
\[
(S_{\al W.}^p)^* \ R \ S_{ \mt{\al W.}.}^p
=
(S_{\al V.}^p)^* \ R \ S_{ \mt{\al V.}.}^p
\ .
\]
\end{lem}
\begin{beweis}
Using Lemma~\ref{lem_fk} and linearity it is enough to consider
expressions of the form
\[
 R = Z R_1 \cdots R_n \;,\; Z \in \al Z.\;,\; R_i \in
\al L. ( \al H._{V_i}^{r(V_i)} , \al H._{V_i}^{r(V_i)+k(V_i)} )\,,
\]
with $\left\{  V_1 , \ldots , V_n  \right\}\subseteq {supp} (R)$.
Now, by Lemma~\ref{InvVect}~(ii) it follows that $S_W^* R_i S_W = R_i$
for every $W \neq V_i$.
\end{beweis}

Next we give an approximation result for elements
of the form as given in Eq.~(\ref{def_ttilde})
with $T\in\al F.^k$ in terms of series of elements in the
tensor product of suitable Hilbert bimodules
(cf.~Remark~\ref{rem_part1_0}).

\begin{lem}
\label{lem_part1}
Let $T \in$ $\al F.^k$, $k \in \Z_0^\infty$.
Then, for every $\varepsilon > 0$ there exist $n(\varepsilon)$, $p(\varepsilon) \in \N$,
finite sets $\al V._1 , \ldots , \al V._{n(\varepsilon)} \subseteq\gfin$,
and $\Phi_l \in \ot H._{\al V._l}^k$ satisfying the following property:
defining the finite subset $\al W.:=\mathop{\bigcup}\limits_{l=1}^{n(\varepsilon)} \al V._l
\subset\gfin$ we have
\[
\left\| \
\widetilde T_{\al W.,p} -
\sum_{l=1}^{n(\varepsilon)} \Phi_l
\ \right\|
< \varepsilon\;,\quad \mr for~all.~p> p(\varepsilon)
\ \ .
\]
\end{lem}
\begin{beweis}
Recall that by Lemma \ref{lem_fk}
every element $T \in \al F.^k$, $k \in \Z_0^\infty$,
can be written as
\[
T = \sum_{l=1}^\infty R_l
\ \ , \ \
\]
where $R_l \in$
$\al L. ( \ot H._{\mt{\al V._l}.}^r , \ot H._{\mt{\al V._l}.}^{r+k} )$,
$r \in \N_0^\infty$, and $\al V._l \subseteq \gfin$ is a finite set.
Take $n(\varepsilon)$ such that
\[
\left\| T - \sum_{l=1}^{n(\varepsilon)} R_l \right\| <\varepsilon\;.
\]
For $\al W.:=\mathop{\bigcup}\limits_l^{n(\varepsilon)} \al V._l
\subset\gfin$,
$p \in \N$, we have
\[
\widetilde T_{\al W.,p} = \sum_{l=1}^\infty
                                      (S_{\al W.}^p)^* \ R_l \ S_{\al W.}^p
                        = \sum_{l=1}^\infty
                                      (S_{\al V._l}^p)^* \ R_l \ S_{\al V._l}^p\; ,
\]
where for the last equation we have used Lemma~\ref{lem_rv}.
Since
\[
R_l =
\sum_j^{\mbox{\tiny fin}} \prod_{ V \in \mt{\al V._l}.}
R_{V,j}
\ \ ,
\]
with $R_{V,j} \in$ $\al L. ( \ot H._\mt V.^\mt {r(V)}. , \ot H._\mt V.^\mt {r(V)+k(V)}.  )$
we can concentrate our analysis on the single factors appearing in 
the preceding
product. From the proof of Proposition~3.5 in \cite{Doplicher98} we have that for
every $V \in \al V._l$, $l = 1 , \ldots , n(\varepsilon)$,
there exist $p(l) \in \N$ such that, for
all $p>p(l)$,
\[
(S_V^p)^* \ R_{V,j} \ S_V^p  = \varphi_{V,j}\;,\quad
 \mr for~some.\quad
 \varphi_{V,j} \in \ot H._V^{k(V)}\,.
\]
This implies that for every $p > p(\varepsilon) :=$
$\max \left\{ p(l) \right\}_{l=1}^{n(\varepsilon)}$,
\[
(S_{\al W.}^p)^* \ R_l \ S_{\al W.}^p
=
(S_{\al V._l}^p)^* \ R_l \ S_{\al V._l}^p
=
\Phi_l
:=
\sum_j^{\mbox{\tiny fin}} \prod_{ V \in \mt{\al V._l}.}
\varphi_{V,j}
\ ,
\]
where the r.h.s. of the preceding equality belongs to
$\ot H._{\mt{\al V._l}.}^k$. Finally,
\[
\left\|
 \widetilde T_{\al W.,p} - \sum_{l=1}^{n(\varepsilon)} \Phi_l
\right\| =
\left\|
  (S_{\al W.}^p)^* \left(
               T - \sum_{l=1}^{n(\varepsilon)} R_l
              \right)
   S_{\al W.}^p
\right\| \leq
\left\|
T - \sum_{l=1}^{n(\varepsilon)} R_l
\right\| <
\varepsilon
\]
and the proof is concluded.
\end{beweis}

\begin{pro}\label{Minim}
The {\it C*}-dynamical system $( \al F. , \al G. , \beta )$
with fixed-point algebra $\al A.$ constructed
in this section is minimal, i.e.
\[
\al A.' \cap \al F. =\al Z.\,.
\]
Moreover, $\al Z.$ coincides with the center of $\al A.$, i.e.
$\al Z.=\al A.' \cap \al A.$.
\end{pro}
\begin{beweis}
Our proof will be divided in two steps. We will apply techniques from
\cite[Propositions~3.4 and 3.5]{Doplicher98}.

\paragraph {i) The inclusions
$\al Z. \subseteq \al A.' \cap \al A.\subseteq \al A.'\cap\al F.$:}
Since $\al A.$ is stable w.r.t. the $\T^\infty$-action given by $\delta$,
it suffices to verify that $Z \in \al Z.$
commutes with the elements of the spectral subspaces
$\al A.^k := \al F.^k \cap \al A.$. If $T \in \al A.^k$, then
there is a sequence $\left\{ T_l \right\}_l$ of
elements of $^0 \al F.^k$ approximating $T$.
By Lemma~\ref{lem_fk} we can write
$T_l =\sum^{\mbox{\tiny fin}}_{\mt{\al V.}.} Z_{\mt{\al V.}.}^l T_{\mt{\al V.}.}^l \ $,
with $T_{\mt{\al V.}.}^l \in$
$\al L.( \al H._{\mt{\al V.}.}^l  , \al H._{\mt{\al V.}.}^{l+k} )$.
Applying the mean $\ot m._{\al G.}$ over the group action (cf.~Eq.~(\ref{group-mean})),
we conclude that $T$ is approximated by the sequence
\[
\ot m._{\al G.}(T_l)
=
\sum^{\mbox{\tiny fin}}_{\mt{\al V.}.} Z_{\mt{\al V.}.}^l
                      \
                      \ot m._{\al G.} (T_{\mt{\al V.}.}^l) \ .
\]
Note that each $\ot m._{\al G.} (T_{\mt{\al V.}.}^l)$ is
a $\al G.$-invariant element of
$\al L.( \al H._{\mt{\al V.}.}^l  ,
   \al H._{\mt{\al V.}.}^{l+k} )$.
This implies that the set $^0 \al A.^k := \ ^0 \al F.^k \cap \al A.$
is dense in $\al A.^k$. Therefore it is enough to show commutativity
w.r.t.~$^0 \al A.^k$.
Let $T =\sum^{\mbox{\tiny fin}}_{\al V.} Z_{\al V.} T_{\al V.}$
be a generic element of $^0 \al A.^k$.
If $Z \in \al Z.$,
then the same argument as in Proposition~\ref{teo_hrs} implies
$T_{\al V.} Z = Z T_{\al V.}$, so that
\[
TZ =
\sum^{\mbox{\tiny fin}}_{\al V.} Z_{\al V.} T_{\al V.} Z =
\sum^{\mbox{\tiny fin}}_{\al V.} Z_{\al V.} Z T_{\al V.} =
ZT \ .
\]
The preceding equation and the definition of the
group action in Eq.~(\ref{eq_g}) proves the inclusion
$\al Z. \subseteq \al A.' \cap \al A. \subseteq \al A.'\cap\al F.$.

\paragraph{ii) The inclusion $\al A.'\cap\al F. \subseteq \al Z.$:}
From Lemma~\ref{StableRel} we can reduce our analysis to the spectral
subspaces of the relative commutant $\al C.:=\al A.'\cap\al F.$.
Let $C\in\al C.^k:=\al F.^k\cap\al C.$, $k\in\Z^\infty_0$.
For each finite $\al W. \subseteq \gfin$, we consider the isometry
$S_{ \mt{\al W.}.}$ defined as in Eq.~(\ref{eq_SW}). By Lemma~\ref{InvVect}~(i),
$S_{ \mt{\al W.}.}$ belongs to $\al A.$, thus for every $p \in \N$ we have
\[
 \widetilde C_{\al W.,p}
 := (S_{\al W.}^p)^* \ C \ S_{ \mt{\al W.}.}^p
 = C \,.
\]
From the approximation result in Lemma~\ref{lem_part1}
we conclude that
\[
C = \sum_{l=1}^\infty \Phi_l \;,\quad\mr where.\quad
              \Phi_l \in\ot H._l^k  := \ot H._{\al V._l}^k\,.
\]
If $k(V) = 0$ for every $V \in \gfin$, then every
$\Phi_l$ belongs to $\al Z.$, and we conclude
that $\al C.^0 \subseteq$ $\al Z.$.
Finally we show that $\al C.^k=0$
for all $k\not=0$. For this purpose, take
$C \in \al C.^k$, with $k(V) \neq 0$ for some $V \in \gfin$ and
consider the following equivalence relation for the indices:
\[
l \sim m
\ :\Leftrightarrow \
k |_{\al V._l} = k |_{\al V._m}
\ \ .
\]
Note that $l \sim m$ if and only if $\Phi_l \Phi_m^* \in$
$\ot H._l^k (\ot H._m^k)^* \subset$ $\al F.^0$.
Let us consider
\[
CC^* =
\sum_{l,m=1}^\infty \Phi_l \Phi^*_m \ \ ;
\]
since $CC^* \in$ $\al C.^0$, by applying the
projection onto the zero spectral component
(cf.~Eq.~(\ref{eq_m0})) we obtain
\begin{equation}
\label{eq_summq}
CC^* = m_0 (CC^*) =
\sum_{l,m} m_0 (\Phi_l \Phi_m^*) =
\sum_{l \sim m} \Phi_l \Phi_m^*
\ .
\end{equation}
We now apply to the sum (\ref{eq_summq}) the argument used for
generic elements of $\al C.^0$, and conclude that $\Phi_l \Phi_m^* \in\al Z.$
for every $l \sim m$. In particular, we obtain
$\Phi_l \Phi_l^* \in\al Z.$ for every $l \in \N$.
Now, since every $\ot H._\mt V.$ is a nonsingular bimodule
(see Definition~\ref{def_nonsing} and Theorem~\ref{thm_main}), the same is true
for $\ot H._l^k$ (see \cite[p.~273]{Doplicher98}). This implies that
$\Phi_l = 0$ for every $l \in \N$.
We conclude that if $C \in \al C.^k$, $k\not=0$, then $C = 0$
and therefore $\al C. =$ $\al C.^0 =\al Z.$.
\end{beweis}

\begin{teo}\label{FinalEx}
Let $\al G.$ be a compact group, $\al Z.$ a unital Abelian {\it C*}-algebra
and $\alpha\colon\ot C.(\al G.) \rightarrow\aut \al Z.$ a fixed
chain group action. Given the set of finite-dimensional
representations $\gfin$ as in Definition~\ref{def_gfin}, the
{\it C*}-dynamical system $( \al F. , \al G. , \beta )$
with fixed-point algebra $\al A.$ constructed
in this section is minimal (i.e.~$\al A.'\cap\al F.=\al Z.$)
and regular. Moreover, $\al Z.$ coincides with the center of $\al A.$, i.e.
$\al Z.=\al A.' \cap \al A.$, and
every algebraic Hilbert space $\al H._\mt V.$,
$V \in \gfin$ (cf.~Definition~\ref{def_gfin}),
is contained in $\al F.$ with support $\1$.
\end{teo}
\begin{beweis}
In Subsection~\ref{5.1} and \ref{5.2} we have specified the construction
of a {\it C*}-dynamical system $( \al F. , \al G. , \beta )$ where
every algebraic Hilbert space $\al H._\mt V.$, $V \in \gfin$,
is contained in $\al F.$ with support $\1$ (recall the relations
(\ref{31.1})-(\ref{31.4})). The minimality property has been shown
in Proposition~\ref{Minim}.

To show that $( \al F. , \al G. , \beta )$ is regular
recall Proposition~\ref{prop0} and Definition~\ref{ReguCondi}.
By construction, the free $\al G.$-invariant bimodules
$\ot H._\mt V.$ are generated by the corresponding
$\al G.$-invariant algebraic Hilbert spaces $\al H._\mt V.$.
From relation (\ref{31.3}) it is clear that the
assignment $\ot H._\mt V.\mapsto\al H._\mt V.$ is compatible with
products, hence $( \al F. , \al G. , \beta )$ is regular.
\end{beweis}

\begin{eje}\label{SU2}
\textbf{(Hilbert C*-systems for SU(2))}\\
We will apply here the preceding theorem to the group $\al G.=\mr SU.(2)$
and a given Abelian {\it C*}-algebra $\al Z.$.
Note that the all nontrivial irreducible representations of $\al G.$
satisfy the properties (i) and (ii) in Definition~\ref{def_gfin}. Indeed,
any irreducible representation
$V^{(l)}$, $l\in\{\ms \frac12.,1,\ms \frac32.,\dots\}$, has dimension
$2l+1\geq 2$ and the decomposition of the tensor product
$V^{(l)}\otimes V^{(l)}$ always contains the trivial representation
$\iota$ as a subrepresentation (cf.~Example~\ref{SU2_v1}). Therefore
we can use $\wg$ to construct the {\it C*}-algebra $\al F.$ as in
Subsection~5.1. Note that since the chain group of $\mr SU.(2)$
is isomorphic to $\Z_2$, there are only two choices of chain group
homomorphism $\alpha\colon\Z_2\to\al Z.$.

Applying Theorem~\ref{FinalEx} we obtain a minimal Hilbert {\it C*}-system
$(\al F.,\mr SU.(2))$
(i.e.~a minimal and regular {\it C*}-dynamical system where all irreducible
representations of the group are realized on algebraic Hilbert spaces with
support $\1$). Note finally that the trivial representation
is realized on the space $\C\1\subset\al Z.$.
\end{eje}

\section{Appendix: Tensor categories of Hilbert bimodules}

Let $\ot R.$ be a {\it C*}-algebra, $\ot H.$ a Hilbert $\ot R.$-bimodule.
We consider the category with objects the internal tensor powers $\ot H.^r$,
$r \in \N_0$, and arrows the spaces $\al L. (\ot H.^r , \ot H.^s)$. It is
well-known that such a category is not a tensor category: in fact, in
general it is not possible to define in a consistent way the tensor
product of right $\ot R.$-module operators
$T \in \al L. (\ot H.^r , \ot H.^s)$,
$T' \in \al L. (\ot H.^{r'} , \ot H.^{s'})$
(see \cite[VI.13.5]{bBlackadar98}).
Nevertheless, the tensor product
$T \otimes T' \in \al L. (\ot H.^{r+r'} , \ot H.^{s+s'})$
makes sense if $T$, $T'$ are $\ot R.${\em -bimodule operators}, i.e.
if they commute with the left action: $T A \psi = A T \psi$,
$A \in \ot R.$, $\psi \in \ot H.^r$.
Thus, the above-mentioned problem of tensoring right
$\ot R.$-module operators vanishes, if we consider
an $\ot R.$-bimodule $\ot H.$ such that $A \psi = \psi A$,
$\psi \in \ot H.$, $A \in \ot R.$. Such a bimodule is called
{\em symmetric}. In the context of diagonal bimodules considered in
this paper symmetry is guaranteed if we have a trivial chain group
action (cf.~Eq.~(\ref{31.3})).

Let $\al G.$ be a compact group such that $\ot H.$ is a
$\al G.$-Hilbert $\ot R.$-bimodule (in the sense of
\cite[VIII.20.1]{bBlackadar98}, \cite[\S 2]{Kasparov80}).
Suppose that every $g \in \al G.$ acts on $\ot H.$ as an
$\ot R.$-bimodule operator. Then, for every $r \in \N_0$ it makes
sense to consider the unitary
$g_r := g \otimes \ldots \otimes g \in \al L. ( \ot H.^r , \ot H.^r )$,
so that there is an action
$\al G. \ni g \rightarrow g_r \in \al L.(\ot H.^r , \ot H.^r )$,
and $\ot H.^r$ is a $\al G.$-Hilbert $\ot R.$-bimodule.
We denote by $\al L. (\ot H.^r , \ot H.^s ; \al G. )$
the spaces of $\al G.$-equivariant operators
$T \in \al L. (\ot H.^r , \ot H.^s)$ such that
$T g_r \psi = g_s T \psi$, $g \in \al G.$, $\psi \in \ot H.^r$.
We denote by ${\bf{tens}} ( \al G. , \ot H. )$ the {\it C*}-category
with objects $\ot H.^r$, $r \in \N_0$, and arrows
$\al L. (\ot H.^r , \ot H.^s ; \al G. )$. For $r=0$, we define
$\ot H.^0 := \ot R.$, so that
\begin{equation}
\label{eq_ag}
\al L. (\ot R. , \ot R.) =
\al L. (\ot R. , \ot R. ; \al G. ) =
\ot R. \ .
\end{equation}
Note that $Z \in \al L. (\ot R. , \ot R.)$
 is an $\ot R.$-bimodule operator if and only if $Z$
 belongs to the center of $\ot R.$.
It is an interesting question to ask whether there exist
non-symmetric, $\al G.$-Hilbert $\ot R.$-bimodules
$\ot H.$ such that ${\bf{tens}} ( \al G. , \ot H. )$ is a
{\em tensor} {\it C*}-category (if $\ot H.$ is symmetric,
it is clear that the above question is trivial).
From the above considerations,
it follows that in order to get a tensor structure on
${\bf{tens}} ( \al G. , \ot H. )$, the following two conditions
are needed:
\begin{enumerate}
\item  if $T \in \al L. (\ot H.^r , \ot H.^s ; \al G. )$, then
       $T$ must be an $\ot R.$-bimodule operator. In fact, in such
       a case it makes sense to consider the tensor product
       $T \otimes T'$, for every
       $T' \in \al L. ( \ot H.^{r'} , \ot H.^{s'} ; \al G.)$;
\item  as a consequence of the preceding point, from
       (\ref{eq_ag}) and subsequent remarks we conclude that $\ot R.$
       has to be Abelian.
\end{enumerate}

A $\al G.$-action satisfying the above properties
is called {\em tensor action}, in the terminology of
\cite[Definition~5.3]{Vasselli05}. The next result
is just a reformulation of Theorem~\ref{thm_main}, and
shows that the class of Hilbert bimodules
admitting a tensor action is quite rich.

\begin{pro}
\label{prop_ex_bimod}
Let $\al G.$ be a compact group, $\al Z. \neq \C\1$ a nontrivial unital Abelian
{\it C*}-algebra with a fixed (nontrivial)
action $\alpha  \colon \ \ot C.(\al G.) \to \mr Aut.(\al Z.)$. Then, there exists
at least a (nonsymmetric) Hilbert $\al Z.$-bimodule $\ot H.$
carrying a tensor $\al G.$-action by unitary $\al Z.$-bimodule operators.
\end{pro}

\paragraph{Acknowledgments}
We are grateful to the DFG-Graduiertenkolleg
``Hierarchie und Symmetrie in mathematischen Modellen''
for supporting a visit of E.V.~to the RWTH-Aachen. It is also a pleasure to
thank Ricardo P\'erez Marco for useful comments on the manuscript.



\end{document}